%% file: HB_WQ_paper_draft1.tex
\documentclass[preprint,1p,fleqn]{elsarticle}
\usepackage{amsfonts,amsmath,amsthm,amssymb}
\usepackage{verbatim}
\usepackage{color}
\usepackage{moreverb}
\usepackage{multirow}

\usepackage[ruled,longend]{algorithm2e}

\usepackage{subfigure}
\usepackage{caption}
\usepackage{wrapfig}
\usepackage{epsfig}
\usepackage{subfig}
\usepackage{mathtools}

\usepackage[normalem]{ulem}

\usepackage{siunitx}
\usepackage{pgfplotstable,booktabs,array,lscape}
\usepackage{pgfplots}
\pgfplotsset{compat=1.16}
\usetikzlibrary{external}
\tikzset{external/force remake}

\usepackage{filecontents}

\pgfplotsset{select coords between index/.style 2 args={
    x filter/.code={
        \ifnum\coordindex<#1\fi
        \ifnum\coordindex>#2\fi
}}}

\def\RR{\mathbb{R}}

\def\bx{\boldsymbol{x}}

\def\bi{\boldsymbol{i}}

\def\m{m}
\def\L{L}
\def\admissibility{r}

\def\c{c}

\DeclareMathOperator\suppoperator{supp}
\newcommand{\supp}[1]{{\suppoperator (#1)}}

\definecolor{color1}{rgb}{0.6,0,0.6}
\definecolor{color2}{rgb}{0.8,0,0}
\definecolor{color3}{rgb}{0.0,0.5,0.0}
\definecolor{color4}{rgb}{0,0,0.7}
\definecolor{color5}{rgb}{0,0.5,0.5}
\definecolor{dgreen}{rgb}{0,0.6,0.0}

\newtheorem{prop}{Proposition}
\newtheorem{remark}{Remark}

\newcommand{\C}[1]{\ensuremath{\mathcal{#1}}}

\newcommand{\bh}{\ensuremath{ {b}}}

\newcommand{\fnTid}[1]{\ensuremath{\boldsymbol{#1}}}
\newcommand{\fnGid}[1]{\ensuremath{#1}}

\newcommand{\activeFnTids}[1]{\ensuremath{
\ensuremath{\C{I}^{#1}_{\mathcal{H}}}}}

\newcommand{\ptTid}[1]{\ensuremath{\mathbf{#1}}}
\newcommand{\ptGid}[1]{\ensuremath{\mathrm{#1}}}

\newcommand{\bspline}[1]{\ensuremath{{\C{B}}^{#1}}}

\newcommand{\bsplineUnivar}[2]{\ensuremath{\bspline{#1}_{#2}}}
\newcommand{\bsplineFnTids}[1]{\ensuremath{\C{I}^{#1}_{\mathcal{B}}}}

\newcommand{\nestedDomain}[1]{\ensuremath{ \Omega^{#1}}}

\newcommand{\hr}[0]{\ensuremath{ \mathcal{H}}}

\newcommand{\hrLvl}[1]{\ensuremath{ \hr^{#1}}}

\newcommand{\ndofsUnivar}[2]{\ensuremath{N^{#1}_{\bspline{}_{#2}}}}

\newcommand{\ndofsBspline}[1]{\ensuremath{N_{\bspline{}}^{#1}}}

\newcommand{\ndofsHr}[0]{\ensuremath{N_{\hr}}}

\newcommand{\ptsLvlTP}[1]{\ensuremath{\mathcal{Q}^{#1}}}

\newcommand{\ptsTPTids}[1]{\ensuremath{\mathcal{I}^{#1}_{\mathcal{Q}}}}

\newcommand{\ptbUnivar}[2]{\ensuremath{\bar{x}^{#1}_{#2,q_{#2}}}}

\newcommand{\ptsUnivar}[2]{\ensuremath{\ptsLvlTP{#1}_{#2}}}

\newcommand{\nptsUnivar}[2]{\ensuremath{R^{#1}_{#2}}}

\newcommand{\nptsLvlTP}[1]{\ensuremath{R^{#1}}}

\newcommand{\fnUnivar}[2]{\ensuremath{ \bh^{#1}_{#2} } }

\newcommand{\fnLvl}[2]{\ensuremath{ 
b^{#1}_{ \fnTid{#2} } } }

\newcommand{\quadrule}[2]{\ensuremath{\mathfrak{Q}^{#1}_{\fnTid{#2}}}}

\newcommand{\dofsHr}[0]{\ensuremath{\mathcal{I}_{\hr}}}

\newcommand{\wtMultivar}[2]{\ensuremath{w^{#1}_{\fnTid{#2},\ptTid{q}}}}
\newcommand{\wtbMultivar}[2]{\ensuremath{\bar{w}^{#1}_{\fnTid{#2},\ptTid{q}}}}
\newcommand{\ptbMultivar}[1]{\ensuremath{\bar{\bx}^{#1}_{\ptTid{q}}}}

\SetKwInOut{Input}{Input}
\SetKwInOut{Output}{Output}

\SetKwFunction{classify}{classify\_basis\_functions}
\SetKwFunction{activepts}{compute\_active\_points}
\SetKwFunction{weightsunivar}{compute\_1Dweights}
\SetKwFunction{activeunivaridinblock}{compute\_active\_1DBsplines}

\begin{document}

\begin{frontmatter}

\title{Weighted quadrature  for hierarchical B-splines}


\author[label1]{Carlotta Giannelli}
\ead{carlotta.giannelli@unifi.it}
%
\author[label2]{Tadej Kandu\v{c}}
\ead{tadej.kanduc@fmf.uni-lj.si}
%
\author[label3]{Massimiliano Martinelli}
\ead{martinelli@imati.cnr.it}
%
\author[label4,label3]{Giancarlo Sangalli}
\ead{giancarlo.sangalli@unipv.it}
%
\author[label3]{Mattia Tani}
\ead{mattia.tani@imati.cnr.it}
%

\address[label1]{Dipartimento di Matematica e Informatica ``U. Dini'', Universit\`a degli Studi di Firenze, Italy }
\address[label2]{Faculty of Mathematics and Physics, University of Ljubljana, Ljubljana, Slovenia}
\address[label3]{Istituto di Matematica Applicata e Tecnologie Informatiche ``E. Magenes'' -- CNR, Pavia, Italy}
\address[label4]{Dipartimento di Matematica ``F. Casorati'', Universit\`a degli Studi di Pavia, Italy}

\begin{abstract}
{We present weighted quadrature for hierarchical B-splines to address the fast formation of system matrices arising from adaptive isogeometric Galerkin methods with suitably graded hierarchical meshes.
By exploiting a local tensor-product structure, we extend the
construction of weighted rules from the tensor-product to the
hierarchical spline setting. The proposed algorithm has a
computational cost proportional  to the number of degrees of freedom and advantageous properties with increasing spline degree. To illustrate the performance of the method and confirm the theoretical estimates, a selection of 2D and 3D numerical tests is provided.}
\end{abstract}

\end{frontmatter}

\section{Introduction}\label{sec:intro}
\input{HB_WQ_Intro.tex}

\section{Preliminaries}\label{sec:pre}
\subsection{Hierarchical B-splines}\label{sec:HB}
\input{HBsplines.tex}
\subsection{WQ for tensor-product splines}\label{sec:WQR}

\input{HB_WQ_TP.tex}

\section{WQ and mass matrix formation
  for hierarchical B-splines}\label{sec:HWQR}
\subsection{Definition of WQ  for hierarchical B-splines}

\input{HB_WQ_Definition.tex}

\subsection{Preprocessing: computing the quadrature points and weights}\label{sec:preprocessing}
\input{HB_WQ_Preprocessing.tex}

\subsection{Matrix formation: algorithm}\label{sec:formation_algo}

\input{HB_WQ_Assembly.tex}

\section{Computational cost} \label{sec:formation_cost}

\input{Cost}

\section{Numerical tests}\label{sec:tests}
\input{HB_WQ_NumericalTests.tex}

\section{Closure}\label{sec:closure}

A fast matrix formation technique for adaptive isogeometric Galerkin methods with multivariate hierarchical B-splines was presented by focusing on the efficient design of weighted quadrature rules. The theoretical estimates of the computational cost suitably exploit the limited number of basis functions which are non-zero on any element of an admissible hierarchical mesh. A selection of numerical examples confirm that the results obtained with the hierarchical weighted approach compare favorably with respect to standard Gaussian quadrature rules, specially in the three-dimensional case. 
Interesting topics for future research include for example  the
combination of the proposed algorithm with matrix-free methods
\cite{st2018} as well as the extension to the case of truncated
hierarchical B-splines \cite{giannelli2012,giannelli2016}, and the application to PDE problems of
applicative interest.

\bibliography{HB_WQ_Ref}

\end{document}

%% file: HB_WQ_Intro.tex
{ 
Local and adaptive mesh refinement methods in isogeometric analysis have gained a notable attention in the last years and their mathematical theory has recently been established, see \cite{bggpv2021} and references therein. One of the more prominent tool in this context is provided by hierarchical B-spline constructions \cite{vuong2011,giannelli2012,giannelli2016}}. The attractive advantage of the hierarchical spline model comes from a good balance between sound theoretical foundations, flexibility, and ease of implementation.  A local refinement step is governed by simple conditions that activate/deactivate basis functions from hierarchically nested sequence of spline spaces. The use of the hierarchical approach in isogeometric analysis was originally proposed in \cite{vuong2011} and subsequently investigated in different directions, which range from the theory of adaptive methods \cite{bg2016,bg2017,ghp2017} to engineering applications, see e.g., \cite{hennig2018,cgrv2019,kuru2013} and references therein.


{ 

The efficient formation of matrices in isogeometric Galerkin methods is a topic of active research. In this paper we focus on the weighted quadrature (WQ) approach,  introduced in \cite{BWR2017}. Other recent results and methods in this area are  integration by interpolation and look-up \cite{pan2020fast,pan2021efficient}, multiscale quadrature \cite{hirschler2021fast},  sum factorization \cite{antolin2015efficient,bressan2019sum}, the surrogate matrix method \cite{drzisga2020surrogate},  reduced integration at superconvergent points \cite{fahrendorf2018reduced} and, beyond quadrature,  the use of low-rank approximation \cite{mantzaflaris2017low} or GPUs \cite{karatarakis2014gpu}.

The aim of WQ is  to reduce the number of
quadrature points that are needed to accurately compute integrals involving products of  B-spline  basis
functions. In combination with sum-factorization and other
implementational techniques, it reduces significantly the cost of
formation of isogeometric matrices.
The idea of WQ is that the test function plays the role of weight function in the
integration, and therefore the quadrature weights depends on the test function. 
 The advantage of this construction is that the number of exactness conditions to be
imposed is less than for Gauss quadrature,  generalized Gauss
quadrature 
\cite{hughes2010efficient,auricchio2012simple,bartovn2017gauss,bartovn2016optimal} or
reduced quadrature \cite{hiemstra2017optimal,schillinger2014reduced,adam2015selective}.
Therefore, WQ requires less quadrature points, which   mildly
depend on the spline degree.}

 { In this paper we extend WQ to hierarchical B-splines with  maximum regularity. Since the construction of the hierarchical basis is simply based on a suitable selection of standard B-splines at different levels of details, we can define hierarchical  WQ as a linear combination of standard WQ on
different tensor-product levels.
The proposed algorithm has a computational cost proportional to the number of degrees of freedom and advantageous properties with increasing spline degree. To illustrate the performance of the method and confirm the theoretical estimates, a selection of 2D and 3D numerical tests is provided. For the sake of simplicity, we discuss the case of the mass matrix and $L^2$-projection. However, dealing with other matrices and with PDE problems is conceptually the same, see \cite{BWR2017}.

The structure of the paper is as follows. Preliminaries on hierarchical B-splines and weighted quadrature are recalled in Section~\ref{sec:pre}. The WQ and its use in the mass matrix formation for hierarchical B-splines is then presented in Section~\ref{sec:HWQR}, while its computational cost is studied in Section \ref{sec:formation_cost}. Section~\ref{sec:tests} illustrates the numerical experiments and,  finally, Section~\ref{sec:closure} concludes the paper.}

%% file: HBsplines.tex

We consider a nested sequence of $L+1$ multivariate tensor-product spline spaces $V^\ell$, for $\ell=0,\ldots,\L,$ of fixed degree $p$ in any coordinate direction defined on a bounded closed hyper-rectangle $\Omega \subset \mathbb{R}^d$. By focusing on dyadic mesh refinement, we assume the spline spaces defined on a sequence of suitably refined knot vectors so that
$V^\ell \subset V^{\ell+1}$, for $\ell=0,\ldots,L-1.$ It should be noted however that the hierarchical B-spline model can be considered also in connection with more general (non-uniform) mesh refinement rules, where each mesh element is subdivided in an arbitrary number of children elements. Moreover, not only $h$-refinement but also $p$-refinement can be combined with the construction of the spline hierarchy as long as the spaces remain nested between each pair $(\ell,\ell+1)$ of consecutive levels, for $\ell =0,\ldots,L-1$.
The considered choice is dyadic (uniform) refinement and fixed spline degree at all levels, which is the standard setting for adaptive isogeometric methods, based on hierarchical B-splines, and a suitable compromise between accuracy and efficiency for related application algorithms. Note that the design and development of fast assembly and efficient numerical integration rules tailored on hierarchical B-spline constructions are key ingredients for the subsequent development of more flexible adaptive approximation schemes.

In direction $k$ of the domain $\Omega$, the level $\ell$ basis  $\bsplineUnivar{\ell}{k}$ consists of $\ndofsUnivar{\ell}{k}$ univariate B-splines $\fnUnivar{\ell}{k,i_k}$,
\begin{equation*}
 \bsplineUnivar{\ell}{k} = \bigl\{ \fnUnivar{\ell}{k,i_k} \; : \; i_k = 1, 2, \dots, \ndofsUnivar{\ell}{k} \bigr\}.
\end{equation*}
%
The multivariate spline space $V^\ell$ on $\Omega$ can be defined as the span of the tensor-product B-spline basis functions $\fnLvl{\ell}{i}$,
\begin{equation*}
 \bspline{\ell} := \bigl\{  \fnLvl{\ell}{i} := \prod_{k=1}^{d} \fnUnivar{\ell}{k,i_{k}} \; : \; \fnUnivar{\ell}{k,i_{k}} \in \bsplineUnivar{\ell}{k} \bigl\} 
\end{equation*}
with respect to the index set
\begin{align*}
 \bsplineFnTids{\ell} 
 & :=  \bigl\{ \fnTid{i} := (i_1,\dots,i_d) \in \mathbb{N}^{d} \; : \; 1 \le i_k \le \ndofsUnivar{\ell}{k}, \, k=1,\dots,d  \bigr\}.
\end{align*}
We denote the rectilinear mesh grid of level $\ell$ by ${\cal M}^\ell$.
The dimension of $V^\ell$ is then simply given by
\begin{align*} \ndofsBspline{\ell} := \dim \bigl( \bspline{\ell} \bigr) = \prod_{k=1}^{d} \dim \bigl( \bsplineUnivar{\ell}{k} \bigr) = \prod_{k=1}^{d} \ndofsUnivar{\ell}{k}\; .
\end{align*}

To localize the refinement regions at different hierarchical levels, we also consider a nested sequence of closed subsets of $\nestedDomain{0} := \Omega$ given by
\begin{align*}
\nestedDomain{0} \supseteq \nestedDomain{1}
\supseteq \ldots \supseteq \nestedDomain{L+1} = \emptyset .
\end{align*}
The hierarchical mesh ${\cal M}$ collects the grid elements $M \in \cal M^\ell$, which are not included in any refined region $\nestedDomain{m}$ of higher level $m>\ell$,
\begin{align*}
{\cal M} := \left\{M\in {\cal M}^\ell \; : \; M\subseteq\nestedDomain{\ell}, M\not\subseteq\nestedDomain{\ell+1},\, \ell=0,\ldots,L\right\}.
\end{align*}

We define the hierarchical B-spline basis \cite{vuong2011} with respect to the hierarchical mesh ${\cal M}$ as
\begin{align*}
\hr({\cal M}) := \left\{
\fnLvl{\ell}{i} \in \bspline{\ell} \; : \; \fnTid{i} \in \activeFnTids{\ell}, \,
\ell=0,\ldots,L
\right\},
\end{align*}
where
\begin{align*}
\activeFnTids{\ell} := \left\{
\fnTid{i} \in \bsplineFnTids{\ell} \; : \; \supp{\fnLvl{\ell}{i}}
\subseteq \Omega^\ell, \supp{\fnLvl{\ell}{i}}\not\subseteq \Omega^{\ell+1}
\right\}.
\end{align*}
The cardinality is denoted by $\ndofsHr$. 
Each basis function $\fnLvl{\ell}{i} \in \hr$ is uniquely identified by its level $\ell$ and by the 
multi-index $\fnTid{i} \in \activeFnTids{\ell}$. 
Hence we can define the set of basis identifiers 
\begin{align} \label{eq:hr_ids}
\dofsHr := \bigcup_{\ell=0}^{L} \{ (\ell,\fnTid{i}) \; : \; \bi \in \activeFnTids{\ell} \} \; .
\end{align}
Hierarchical B-splines are non-negative, linear independent, and allow localized mesh refinement by suitably selecting basis functions with a varying level of resolution. 

In order to limit the interaction between B-splines introduced at very different levels of the spline hierarchy, we consider \emph{admissible} meshes. A hierarchical mesh ${\cal M}$ is admissible of class $r$, with $2\le r< L+1$, if the hierarchical B-splines taking non-zero values on any element $Q\in {\cal M}$ belong to at most $r$ successive levels. We refer to \cite{bg2016,bggpv2021} for more details concerning admissible meshes and their properties, and to \cite{bracco2018b} for the presentation of the refinement algorithms which guarantee the construction of hierarchical mesh configurations with different class of admissibility. Note that the suitably graded meshes generated via these algorithms are characterized by a different (stronger) version of admissibility, which is easier to obtain via automatically-driven refinement rules.


%% file: HB_WQ_TP.tex
The goal of WQ is to reduce the number of
quadrature points  in  the computation of Galerkin integrals  for smooth  B-spline  basis
functions. In combination with the sum-factorization and replacement of the
element-wise assembly loop by a direct function-wise calculation of the matrix
entries, the cost  for the formation of isogeometric Galerkin matrices
goes from  $O(p^{3d}N)$ FLOPS down to $O(p^{d+1}N)$ FLOPS, where $p$ and $N$ denote the degree and
degrees-of-freedom number, respectively.

Consider $\bspline{\ell}$,  the level $\ell$   tensor product
B-spline basis, then we can directly apply  the construction of
\cite{BWR2017} and introduce for each $\fnTid{i} \in
\mathcal{I}_{\mathcal{B}}^{\ell} $, the following WQ 
\begin{align}
  \label{eq:WQ-TP-level_l}
  \bar{\mathfrak{Q}}^\ell_{\fnTid{i}}(v):= \sum_{\ptTid{q}
  \in \ptsTPTids{\ell}} \wtbMultivar{\ell}{i} v \left(
  \ptbMultivar{\ell} \right) \approx \int_{[0,1]^d} v(\bx)  \fnLvl{\ell}{i} (\bx) d\bx 
\end{align}
where $\fnLvl{\ell}{i} \in \bspline{\ell} $, $\wtbMultivar{\ell}{i}$ are the quadrature weights and  $\ptbMultivar{\ell} $   the quadrature
points. For a better distinction between tensor-product and hierarchical objects, we use bars on the top of symbols for quadrature rules, points and weights in the former case.
Note that weights depend on $\fnTid{i}$, that is, on $ \fnLvl{\ell}{i} $,  which plays the role of a weight function for the integration.
The rule  \eqref{eq:WQ-TP-level_l} can be used to approximate the $( \fnTid{i} ,  \fnTid{j})$-entry of 
the mass matrix at level $\ell$, indeed
\begin{equation}\label{eq:wq-mass}
  \int_{[0,1]^d} \c (\bx)  \fnLvl{\ell}{j} (\bx)  \fnLvl{\ell}{i}
  (\bx) d\bx  \approx  \bar{\mathfrak{Q}}^\ell_{\fnTid{i}}(\c  \fnLvl{\ell}{j} ) , 
\end{equation}
where  the given function  $c$ takes into account the determinant of
the Jacobian of the parametrization map. More generally, any  matrix
arising  in a isogeometric Galerkin method can be formed by suitable WQ, see \cite{BWR2017}. 
The accuracy of WQ is related to the  exactness conditions that the rule satisfies. 
In the case of \eqref{eq:WQ-TP-level_l}--\eqref{eq:wq-mass}, a typical request is
\begin{equation}\label{eq:exactness-of-wq-in-1D}
  \bar{\mathfrak{Q}}^\ell_{\fnTid{i}}( \fnLvl{\ell}{j} )  =\int_{[0,1]^d} \fnLvl{\ell}{j} (\bx)  \fnLvl{\ell}{i}
  (\bx) d\bx  
  \qquad \forall  \fnLvl{\ell}{j}  \in \bspline{\ell}. 
\end{equation}

Though not necessary, following \cite{BWR2017}, the quadrature points $\ptbMultivar{\ell}$ 
in \eqref{eq:WQ-TP-level_l} do not depend on $\fnTid{i}$.
The set of $d$-variate quadrature points is defined as the following tensor-product 
\begin{equation} \label{eq:ptsTP}
  \ptsLvlTP{\ell} := \ptsUnivar{\ell}{1} \times \dots \times \ptsUnivar{\ell}{d} \; ,
\end{equation}
where 
$\ptsUnivar{\ell}{k} := \{ \bar x_{k,q_k}^{\ell} \}_{q_{k}=1}^{\nptsUnivar{\ell}{k}}$ 
is the set of univariate
quadrature points in the $k$-th direction and $\nptsUnivar{\ell}{k}$
is the number of quadrature points along that direction. In the case of splines of
of maximum regularity, these points can be selected as the
midpoints and endpoints of the knot spans, with the exception of the
first and last knot spans where $p+1$ uniformly distributed quadrature
points can be selected at each edge of the interval, see \cite{BWR2017}.
We also introduce the tensor product set of multi-indices, associated to $\ptsLvlTP{\ell}$ as 
\begin{equation*} \label{eq:idsTP}
  \ptsTPTids{\ell} := \{ \ptTid{q} = (q_{1},\dots,q_{d}) \in \mathbb{N}^{d} \; : 
  \; 1 \le q_{k} \le  \nptsUnivar{\ell}{k} \, , \, 1 \le k \le d \} \; .
\end{equation*}
Even though the quadrature points are defined globally, only those in the support of $ \fnLvl{\ell}{i} $ are active for $\bar{\mathfrak{Q}}^\ell_{\fnTid{i}}$, thus the active quadrature points of $\bar{\mathfrak{Q}}^\ell_{\fnTid{i}}$  depend on $\fnTid{i}$. 
This is formalized by setting to 0 the weights $\{ \wtbMultivar{\ell}{i}
\}_{\ptTid{q}  \in \ptsTPTids{\ell}}$  corresponding to points that are outside the support of $\fnLvl{\ell}{i}$. The nonzero weights are computed by imposing the univariate local exactness conditions, leading to linear problems, whose solution cost is not prevailing in the overall matrix formation cost.





%% file: HB_WQ_Definition.tex
A weighted quadrature rule, associated to an active basis function $\fnLvl{\ell}{i} \in \hr$, is denoted by $\quadrule{\ell}{i}$. Its quadrature points and weights are jointly indexed with respect to an index set that is denoted by $\ptsTPTids{(\ell,\fnTid{i})}$. 
Namely, let 
\begin{equation}
  \label{eq::pts_def}
  \ptsLvlTP{\left(\ell,\fnTid{i}\right)} \coloneqq \left\{ \bx^{\ell}_{\fnTid{i}, \ptTid{q}} 
  \right\}_{\ptTid{q} \in \ptsTPTids{(\ell,\fnTid{i})}}
\end{equation}  
be a set of quadrature points and the set of the corresponding weights are $\left\{ \wtMultivar{\ell}{i} \right\}_{\ptTid{q} \in \ptsTPTids{(\ell,\fnTid{i})}}$. 
The quadrature rule  $\quadrule{\ell}{i}$ applied to an auxiliary function $v$ has the following form
\begin{equation}
  \label{eq:rule}
  \quadrule{\ell}{i} (v) := \sum_{\ptTid{q} \in \ptsTPTids{(\ell,\fnTid{i})} } 
  \wtMultivar{\ell}{i} v \left( \bx^{\ell}_{\fnTid{i},\ptTid{q}} \right) \approx 
  \int_{[0,1]^d} \fnLvl{\ell}{i} (\bx) v(\bx) d \bx.
\end{equation}
A peculiar feature of this structure is that both the set of quadrature points and the set of quadrature weights depend on the considered test function $\fnLvl{\ell}{i}$. However, as we will see in the following, $\ptsLvlTP{\left(\ell,\fnTid{i}\right)}$ can be conveniently selected as a subset of a global tensor-product grid, which is chosen a priori.

Similarly as in the non-hierarchical case discussed in the previous section, the quadrature rules are characterized by exactness conditions. 
More specifically, we require that the rules are exact for all functions in the spline space, or equivalently that
\begin{align}
  \label{eq:exactness_cond}
  \quadrule{\ell}{i} (\fnLvl{m}{j}) = \sum_{\ptTid{q} \in \ptsTPTids{(\ell,\fnTid{i})} } 
  \wtMultivar{\ell}{i} \; \fnLvl{m}{j}\left( \bx^{\ell}_{\fnTid{i},\ptTid{q}} \right) = 
  \int_{[0,1]^d} \fnLvl{\ell}{i} (\bx) \fnLvl{m}{j}(\bx) d \bx \qquad \forall \; 
  (m,\fnTid{j}) \in \mathcal{I}_{\mathcal{H}}.
\end{align}

For a given pair $\left(\ell,\fnTid{i}\right) \in \mathcal{I}_{\mathcal{H}}$ 
we define $\nu(\ell,\fnTid{i})$ as the finest level of a hierarchical basis function such 
that its support has a nonempty intersection with the support of $\fnLvl{\ell}{i}$, i.e.,
\begin{align} 
  \label{eq:mu} 
  \nu(\ell,\fnTid{i}) \coloneqq \max \left\{m \; : \; \supp{ \fnLvl{\ell}{i} } 
  \cap \supp{ \fnLvl{m}{j}  } \neq \emptyset, \; \forall \fnLvl{m}{j} \in \hr \right \}.
\end{align}
To make the notation lighter, the argument $(\ell,\fnTid{i})$ in $\nu (\ell,\fnTid{i})$ will be sometimes omitted, 
since the dependence on the basis identifier will be clear from the context.

Any active basis function $\fnLvl{m}{j} $ that interacts with $\fnLvl{\ell}{i}$ 
(including itself) can be written as a linear combination of basis functions of 
level $\nu(\ell,\fnTid{i})$, that is 
\begin{align} \label{eq:alpha_coef}
  \fnLvl{m}{j} = \sum_{\fnTid{t} \in \bsplineFnTids{\nu}} \alpha_{\fnTid{j},\fnTid{t}}^{m,\nu} \; 
  \fnLvl{\nu}{t} \qquad \forall \; (m,\fnTid{j}) \in \mathcal{I}_{\mathcal{H}} \quad 
  \text{such that} \quad \supp{ \fnLvl{\ell}{i} } \cap \supp{\fnLvl{m}{j}} \neq \emptyset,
\end{align}
with $\alpha_{\fnTid{j},\fnTid{t}}^{m,\nu} > 0$ iff $\supp{ \fnLvl{\nu}{t} } \subseteq \supp{\fnLvl{m}{j}}$ 
and $\alpha_{\fnTid{j},\fnTid{t}}^{m,\nu} = 0$ otherwise.

\begin{figure}[t!]
\centering
\subfigure[Hierarchical mesh and three basis functions from $\bspline{\ell}$]
{
\begin{tikzpicture}
\node[inner sep=0pt] (IG) at (0,0)
{\includegraphics[trim = 3.5cm 1.5cm 3cm 1.25cm, clip = true, height=5.5cm]{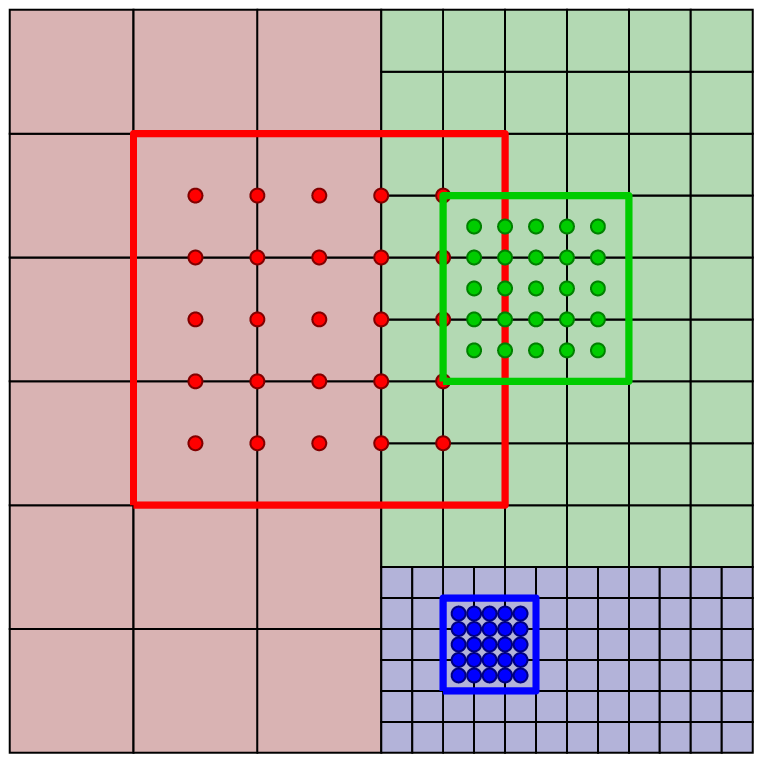}};
{\color{red} \draw[->] (-3.,.25) -- (-1.69,.25);}
{\color{red}\node[text width=0cm] at (-3.5,.25) {$\fnLvl{1}{i}$};}

{\color{dgreen} \draw[->] (-3.,1) -- (.4,1);}
{\color{dgreen} \node[text width=0cm] at (-3.5,1) {$\fnLvl{2}{j}$};}

{\color{blue} \draw[->] (-3.,-1.7) -- (.4,-1.7);}
{\color{blue} \node[text width=0cm] at (-3.5,-1.7) {$\fnLvl{3}{k}$};}

\end{tikzpicture}
}
\subfigure[Hierarchical quadrature points for $\fnLvl{1}{i} \in \hr$]
{
\begin{tikzpicture}
\node[inner sep=0pt] (IG) at (0,0)
{\includegraphics[trim = 3.5cm 1.5cm 3cm 1.25cm, clip = true, height=5.5cm]{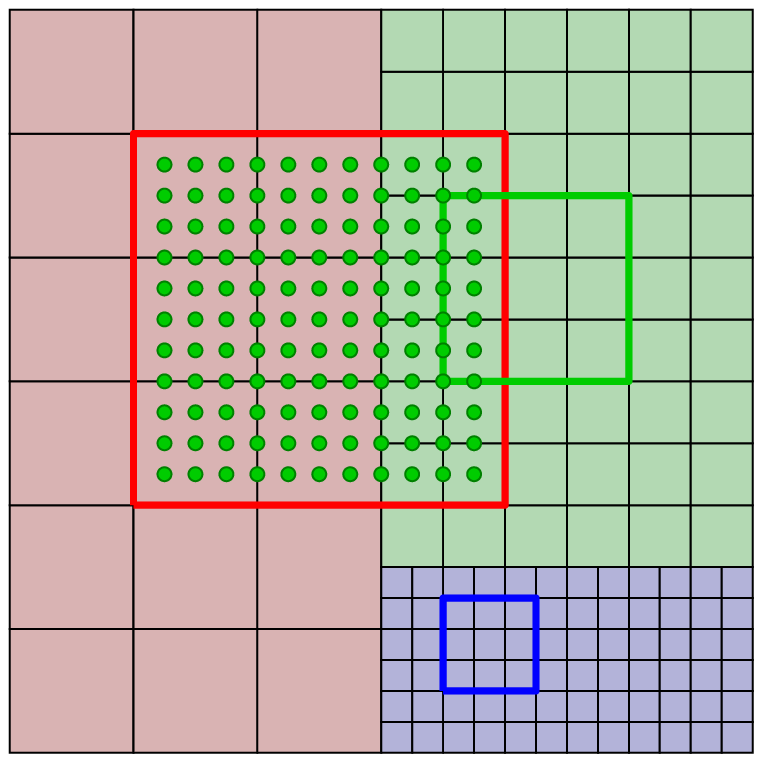}};
{\color{red} \draw[->] (-2.5,.5) -- (-0.4,.5);}
{\color{red} \node[text width=0cm] at (-3.5,.5) {$\ptsLvlTP{\left(1,\fnTid{i}\right)}$};}

\end{tikzpicture}
}
\caption{In (a), the supports and quadrature points for three basis functions from spaces $\bspline{1}$,  $\bspline{2}$ and  $\bspline{3}$ are depicted on a hierarchical mesh, color coded in red, green and blue, respectively. In (b), quadrature points for $\fnLvl{1}{i} \in \hr$ are shown.}
\label{fig:HB_interaction}
\end{figure}

In order to define the quadrature rule $\quadrule{\ell}{i}$ for the hierarchical space $\hr$, 
we rely on the definitions and relations introduced in the previous section for tensor-product spaces. 
In \eqref{eq:rule} we take
\begin{align*}
  \ptsTPTids{(\ell,\fnTid{i})} \coloneqq \left\{ \ptTid{q} \in \ptsTPTids{\nu} \; : \; 
  \ptbMultivar{\nu} \in \ptsLvlTP{\nu} \cap  \supp{\fnLvl{\ell}{i}} \right\} \;
\end{align*}
and let the quadrature points be $ \bx^{\ell}_{\fnTid{i},\ptTid{q}} = \ptbMultivar{\nu}$ 
for every $\ptTid{q} \in \ptsTPTids{(\ell,\fnTid{i})}$, hence (from definition~\eqref{eq::pts_def}) we have
\begin{align}
  \label{eq:hier_pointSet}
  \ptsLvlTP{\left(\ell,\fnTid{i}\right)} = \ptsLvlTP{\nu} \cap \supp{\fnLvl{\ell}{i}} \; .
\end{align}
See Figure~\ref{fig:HB_interaction} for an example of quadrature points for a basis function $\fnLvl{1}{i} \in \hr$. Since it interacts with a level 2 basis function and not with level 3 one, it inherits a local set of level 2 quadrature points \ptbMultivar{2}.

The quadrature weights for $\fnLvl{\ell}{i}$ are simply
\begin{align} 
  \label{eq:hier_weights} 
  \wtMultivar{\ell}{i} =\sum_{\fnTid{s} \in \bsplineFnTids{\nu}} \alpha_{\fnTid{i},\fnTid{s}}^{\ell,\nu} \wtbMultivar{\nu}{s}
\end{align}
for every $\ptTid{q} \in \ptsTPTids{(\ell,\fnTid{i})}$, and the coefficients 
$\alpha_{\fnTid{i},\fnTid{s}}^{\ell,\nu}$ are the expansion coefficients of 
$\fnLvl{\ell}{i}$ on the basis $\bspline{\nu}$ as in \eqref{eq:alpha_coef}. 

In the following proposition we show that this choice for the quadrature rule 
$\quadrule{\ell}{i}$ satisfies the imposed exactness conditions on the hierarchical space.

\begin{prop}
\label{prop:2}

For $(\ell,\fnTid{i}) \in \mathcal{I}_{\mathcal{H}}$ the quadrature rule $\quadrule{\ell}{i}$ satisfies the exactness conditions \eqref{eq:exactness_cond} on the hierarchical space.
\end{prop}

\proof
Let $(m,\fnTid{j}) \in \mathcal{I}_{\mathcal{H}}$. 
If $\supp{ \fnLvl{\ell}{i} } \cap \supp{\fnLvl{m}{j}} = \emptyset$, then the equation in 
\eqref{eq:exactness_cond} is trivially satisfied since the quadrature points $\ptsLvlTP{\left(\ell,\fnTid{i}\right)}$
belong to the support of $\fnLvl{\ell}{i}$.
On the other hand, if $\supp{ \fnLvl{\ell}{i} } \cap \supp{\fnLvl{m}{j}} \neq \emptyset$, then
\begin{align*}
  \quadrule{\ell}{i}\left( \fnLvl{m}{j} \right) 
  & = \sum_{\ptTid{q} \in \ptsTPTids{(\ell,\fnTid{i})} } \wtMultivar{\ell}{i} \fnLvl{m}{j} \left( \bx^{\ell}_{\fnTid{i},\ptTid{q}} \right)
  = \sum_{\ptTid{q} \in \ptsTPTids{(\ell,\fnTid{i})} } \wtMultivar{\ell}{i} \fnLvl{m}{j} \left( \ptbMultivar{\nu} \right)
  \\ & = \sum_{\ptTid{q} \in \ptsTPTids{\nu}} \left( \sum_{\fnTid{s} \in 
  \bsplineFnTids{n}} \alpha_{\fnTid{i},\fnTid{s}}^{\ell,\nu} \wtbMultivar{\nu}{s} \right) 
  \left( \sum_{\fnTid{t} \in \bsplineFnTids{\nu}} \alpha_{\fnTid{j},\fnTid{t}}^{m,\nu} \fnLvl{\nu}{t} 
  \left( \ptbMultivar{\nu} \right) \right)
  \\ & = \sum_{\fnTid{s},\fnTid{t} \in \bsplineFnTids{\nu}} \alpha_{\fnTid{i},\fnTid{s}}^{\ell,\nu} 
  \alpha_{\fnTid{j},\fnTid{t}}^{m,\nu} \sum_{\ptTid{q} \in \ptsTPTids{\nu}} \wtbMultivar{\nu}{s} 
  \fnLvl{\nu}{t} \left( \ptbMultivar{\nu} \right)  
  \\ & = \sum_{\fnTid{s},\fnTid{t} \in \bsplineFnTids{\nu}} \alpha_{\fnTid{i},\fnTid{s}}^{\ell,\nu} 
  \alpha_{\fnTid{j},\fnTid{t}}^{m,\nu} \int_{[0,1]^d} \fnLvl{\nu}{s}(\bx) \fnLvl{\nu}{t}(\bx) d \bx
  \\ & = \int_{[0,1]^d} \left( \sum_{\fnTid{s} \in \bsplineFnTids{\nu}} \alpha_{\fnTid{i},\fnTid{s}}^{\ell,\nu} 
  \fnLvl{\nu}{s}(\bx)\right) \left( \sum_{\fnTid{s} \in \bsplineFnTids{\nu}} \alpha_{\fnTid{j},\fnTid{t}}^{m,\nu} 
  \fnLvl{\nu}{t}(\bx)\right) d \bx = \int_{[0,1]^d}  \fnLvl{\ell}{i}(\bx) \fnLvl{m}{j}(\bx) d \bx.
\end{align*}
In the penultimate row we use the property that for each $\fnTid{s} \in \bsplineFnTids{\nu}$ 
the set of points $\ptsLvlTP{\nu}$ and the set of weights $\left\{ \wtbMultivar{\nu}{s} \right\}_{\ptTid{q} \in \ptsTPTids{\nu}}$ 
satisfy the exactness conditions \eqref{eq:exactness-of-wq-in-1D} on the level $\nu = \nu(\ell,\fnTid{i})$.
\endproof


\begin{remark} \label{remark:1}
Quadrature rule $\quadrule{\ell}{i}$ in \eqref{eq:rule} of level $\ell$ is actually 
a linear combination of quadrature rules defined at level $\nu$. 
Namely, if $\supp{ \fnLvl{\ell}{i} } \cap \supp{v} \neq \emptyset$ we can define
${\cal I} := \{\fnTid{s} \in \bsplineFnTids{\nu} \; : \; \supp{\fnLvl{\nu}{s}} \subseteq \supp{\fnLvl{\ell}{i}} \}$
and derive
\begin{align*} 
  \quadrule{\ell}{i} (v) 
  & = \sum_{\ptTid{q} \in \ptsTPTids{(\ell,\fnTid{i})} } \wtMultivar{\ell}{i} 
  v \left( \bx^{\ell}_{\fnTid{i},\ptTid{q}} \right)
  \\ & = \sum_{\ptTid{q} \in \ptsTPTids{(\ell,\fnTid{i})} } \wtMultivar{\ell}{i} v \left( \ptbMultivar{\nu} \right)  
  \\ & = \sum_{\ptTid{q} \in \ptsTPTids{(\ell,\fnTid{i})} } \sum_{\fnTid{s} \in{\cal I}} 
  \alpha_{\fnTid{i},\fnTid{s}}^{\ell,\nu} \wtbMultivar{\nu}{s} v \left( \ptbMultivar{\nu} \right)
  \\ & = \sum_{\fnTid{s} \in {\cal I}} \alpha_{\fnTid{i},\fnTid{s}}^{\ell,\nu} 
  \sum_{\ptTid{q} \in \ptsTPTids{(\nu,\fnTid{s})} } \wtbMultivar{\nu}{s} v \left( \ptbMultivar{\nu} \right) \label{eq:quad_linear_combination-1}
  \\ & = \sum_{\fnTid{s} \in {\cal I}} \alpha_{\fnTid{i},\fnTid{s}}^{\ell,\nu} 
  \bar{\mathfrak{Q}}^{\nu}_{\fnTid{s}} (v).
\end{align*}
In particular, the quadrature rule $\quadrule{\ell}{i}$ is determined by those rules 
$\bar{\mathfrak{Q}}^{\nu}_{\fnTid{s}}$ whose $\supp{\fnLvl{\nu}{s}}$ are in the support of $\fnLvl{\ell}{i}$.
\end{remark}





%% file: HB_WQ_Preprocessing.tex
 
Since the quadrature weights are not known in advance for every 
possible mesh, degree, level and interaction, they need to be 
computed efficiently in the preprocessing phase, before utilizing 
them in the matrix formation phase. For computational 
efficiency, we fully exploit the tensor-product structure of the 
active basis functions $\fnLvl{\ell}{i}$  and of the quadrature 
points for the full level $\ell$ basis functions. Quadrature 
weights are therefore obtained in two steps. First, we compute the 
univariate quadrature weights of level $\nu$, defined in \eqref{eq:mu}, 
by solving the linear systems arising from the univariate exactness 
conditions analogous to \eqref{eq:exactness-of-wq-in-1D}; 
this is the same as in \cite{BWR2017}. 
Then, the univariate quadrature weights for the WQ associated to an 
active basis function $\fnLvl{\ell}{i}$ are computed as linear 
combination of level $\nu$ quadrature weights, analogously to \eqref{eq:hier_weights}.

The quadrature points and weights for $d$-variate B-splines are stored 
and used as $d$-tuples of the univariate points and univariate weights, 
respectively, in order to be ready for the sum-factorization used in the 
matrix formation.

To avoid redundant computations, all the active basis functions are
clustered with  respect to the value of $\nu$ so that the univariate
routines are engaged only ones for each level,
i.e., we classify the basis functions of $\hr$ with respect to $\nu$ by 
defining the {\emph sets of level $n$ interacting functions}
\begin{equation*}\label{eq:Fm}
  F^n \coloneqq \{ \fnLvl{\ell}{\fnTid{i}} \; : \; (\ell,\fnTid{i}) 
  \in \mathcal{I}_{\mathcal{H}} \, , \, \nu({\ell},\fnTid{i}) = n  \}.  
\end{equation*}
It is trivial to check the following properties: 
\begin{align}
 & \hr = \bigcup_{n \leq L} F^n \label{eq:prop1}, \\
 &  F^m \cap F^n = \emptyset \; \text{if} \; m \ne n \label{eq:prop2}, \\
 & 
 F^n \; \text{only contains functions of} \; \hr \; \text{of levels} \; \le n \;. \label{eq:prop3}
\end{align}
The classification of the active basis function $\fnLvl{\ell}{i}$ with respect 
to the maximum level of interaction $\nu({\ell},\fnTid{i})$ is described in 
Algorithm~\ref{algo:classify}.
\begin{algorithm}
  \DontPrintSemicolon
  \Input{$\hr$}
  \BlankLine
  
  Initialize $F^{n} = \emptyset$ for $n=0,\dots,L$\;
  \ForEach{$ (\ell,\fnTid{i}) \in \dofsHr $}
  {
    $F^{\nu(\ell,\fnTid{i})} = F^{\nu(\ell,\fnTid{i})} \cup \fnLvl{\ell}{i}$\;
  }
  \BlankLine
  
  \Output{$\{ F^{n}\}_{n \leq L}$}  
  \caption{\texttt{classify\_basis\_functions}}
  \label{algo:classify}
\end{algorithm}

Because of \eqref{eq:alpha_coef} and \eqref{eq:prop3}, every function $\fnLvl{\ell}{i} \in F^n$ 
can be written as linear combination of functions $\fnLvl{n}{j}$ from level $n$, 
and a similar formula holds for its quadrature weights (see \eqref{eq:hier_weights}). 
Analogously, the same can be said for the (univariate) components: $ \fnUnivar{\ell}{k,i_k}$ 
can be written as a linear combination of level $n$ functions $\fnUnivar{n}{k,j_k}$, 
and $w_{k,i_k,q_k}^{\ell}$ as a linear combination of $\bar w_{k,j_k,q_k}^n$,
\begin{align}\label{eq:1Dref}
\fnUnivar{\ell}{k,i_k} = \sum_{j_k \in D_{k,i_k}^{\ell,n}} \alpha_{k,{i_k},{j_k}}^{\ell,n} \;  \fnUnivar{n}{k,j_k},
\qquad
w_{k,i_k,q_k}^{\ell} = \sum_{j_k \in D_{k,i_k}^{\ell,n}} \alpha_{k,i_k,j_k}^{\ell,n} \bar w_{k,j_k,q_k}^n
\end{align}
for $k=1,2,\dots,d$, where 
\begin{align*} 
  D_{k,i_k}^{\ell,n} = \{ j_k \in \{ 1,\ldots,\ndofsUnivar{n}{k} \} \; : 
  \; \supp{\fnUnivar{n}{k,j_{k}}} \subseteq \supp{\fnUnivar{\ell}{k,i_{k}}} \}. 
\end{align*}

To switch from the $d$-variate to the univariate setting, we first need to define 
two auxiliary functions $\pi_{k}$ and $\tau_{k}$,
\begin{align*}
    \pi_{k} \fnTid{i} \coloneqq i_{k},\qquad  \tau_{k} \fnTid{i} \coloneqq  (i_{1},\dots,i_{k}),
\end{align*}
acting on a multi-index $\fnTid{i} = (i_{1},\dots,i_{d})$. 
Then, for $k=1,\dots,d$ we introduce the set of indices of univariate B-spline that 
are used to define the functions $\fnLvl{\ell}{i} \in \hrLvl{\ell} \cap F^{n}$:
\begin{equation} \label{eq:bsp_active_id_k}
  G^{\ell,n}_{k} \coloneqq \{ i_{k} = \pi_{k} \fnTid{i} \; : \; \fnLvl{\ell}{i} \in \hrLvl{\ell} \cap F^{n}\}. 
\end{equation}
Finally, we define 
\begin{align}\label{eq:set_Dnk}
  D_k^{n} & \coloneqq \bigcup_{\ell = 0}^n \bigcup_{i_k \in G^{\ell,n}_{k} } D_{k,i_k}^{\ell,n}
\end{align}

The sets $D_k^{n}$ identify the univariate  quadrature weights that
are needed to set up the WQ rule. Precisely, in the next step for each
$j_k \in D_k^{n}$ we compute the nonzero univariate  quadrature weights 
$\{ \bar w_{k,j_k ,q_k}^{n} \}_{q_k=1}^{\nptsUnivar{n}{k}}$ associated
to $\fnUnivar{n}{k,j_k} \in \bsplineUnivar{n}{k}$ by imposing the
univariate exactness conditions analogous to
\eqref{eq:exactness-of-wq-in-1D}, exactly as done in
\cite{BWR2017}. Namely, we impose that
\begin{equation*}\label{eq:zero_w} 
  \bar w_{k,j_k,q_k}^{n} = 0 \qquad \text{if} \; \bar x_{k,q_k}^{n} 
  \notin \supp{\fnUnivar{n}{k,j_k}},
  \end{equation*}
while non-zero quadrature weights are obtained by solving the linear system
\begin{equation} \label{eq:exactness1D}
  \sum_{q_k \in A_{k,j_k}^{n,n}}  \bar w_{k,j_k ,q_k}^{n} \fnUnivar{n}{k,t_k}(\bar x_{k,q_k}^{n}) 
  = \int_0^1 \fnUnivar{n}{k,j_k }(\xi) \, \fnUnivar{n}{k,t_k}(\xi) \, d\xi, 
  \qquad \fnUnivar{n}{k,t_k} \in B_{k,j_k}^{n},
\end{equation}
where
\begin{align}
A_{k,j_k}^{\ell,n} &\coloneqq \{ q_{k} \in \{1,\dots,\nptsUnivar{n}{k}\} \; : \; \ptbUnivar{n}{k} \in \ptsUnivar{n}{k} \cap \supp{\fnUnivar{\ell}{k,j_{k}}} \},
  \label{eq:exactness1D_act_ind} \\
  B_{k,j_k}^{n} &\coloneqq \{ \fnUnivar{n}{k,t_k} \in \bsplineUnivar{n}{k} \; : 
  \; \supp{\fnUnivar{n}{k,j_k}} \cap \supp{\fnUnivar{n}{k,t_k}} \ne \emptyset \} 
  \label{eq:exactness1D_trial_B}
\end{align}    
are the set of indices of quadrature points inside the support of the basis function $\fnUnivar{\ell}{k,j_k}$, and the corresponding interacting trial univariate functions, respectively. 
%
%
%
The construction of the univariate quadrature weights is summarized in Algorithm~\ref{algo:weightsunivar}.

\begin{algorithm} 
  \DontPrintSemicolon
  \Input{$\bsplineUnivar{n}{k}, j_{k}, \ptsUnivar{n}{k}$}
  \BlankLine
  
  compute the indices  $A_{k,j_k}^{n,n}$ from \eqref{eq:exactness1D_act_ind}\; 
  find the interacting trial univariate functions $B_{k,j_k}^{n}$ from \eqref{eq:exactness1D_trial_B}\;
  compute non-zero weights $\{\bar w_{k,j_k ,q_k}^{n}\}_{q_k \in A_{k,j_k}^{n,n}}$ by solving the linear system \eqref{eq:exactness1D}\;
  \BlankLine
  
  \Output{$\{ \bar w_{k,j_k ,q_k}^{n} \}_{q_k \in A_{k,j_k}^{n,n}}$}  
  \caption{\texttt{compute\_1Dweights}}
  \label{algo:weightsunivar}
\end{algorithm}

Up to this point, we have defined the univariate quadrature points and computed the univariate quadrature weights, associated to all the basis functions in $\bsplineUnivar{n}{k}$ that are needed to represent functions in $F^{n}$ as linear combination of functions of level $n$, by using \eqref{eq:alpha_coef}.
Using \eqref{eq:1Dref} we can then compute the level $\ell$ univariate quadrature weights 
{
\begin{equation} \label{eq:linear_comb_w_1D}
\mathcal{W}^{\ell,n}_{k,i_k} \coloneqq \Bigl\{ w_{k,i_k,q_k}^{\ell} = \sum_{j_k \in D_{k,i_k}^{\ell , n}} \alpha_{k,i_k,j_k}^{\ell,n} \bar w_{k,j_k,q_k}^n  \; : \; q_k \in A_{k,i_k}^{\ell,n} \Bigr\} \qquad k=1,\dots,d \; ,
\end{equation}
for each index $i_{k} \in G_{k}^{\ell,n}$.
}

The last preprocessing phase is to define the subset of
$d$-dimensional quadrature points that are contained by the support of
functions in $F^{n}$, that will be used in the matrix formation  phase for the evaluation of the non-tensor product coefficients. The union of support of basis functions in $F^{n}$,
\begin{equation} \label{eq:Psin}
  \Psi^{n} \coloneqq \bigcup_{\fnLvl{\ell}{i} \in F^{n} } \supp{\fnLvl{\ell}{i}},
\end{equation}
is a set of $d$-dimensional boxes in $[0,1]^{d}$ that can be described as a set of mesh cells on level $n$, which in general does not have a tensor-product structure. The $d$-dimensional level $n$ quadrature points are simply defined as
\begin{equation} \label{eq:QPsin}
\ptsLvlTP{n}_{\Psi} \coloneqq \ptsLvlTP{n} \cap \Psi^{n}.
\end{equation}
%
%
%
%
\begin{remark} \label{rem:ptsunion}
  Due to nestedness of quadrature points $\ptsLvlTP{n}$ with respect
  to level $n$, there are configurations in which some points are
  defined in multiple levels. For the sake of efficiency, in our code
  we also store  the union of all $d$-dimensional set of points,
  $\ptsLvlTP{}_{\Psi} \coloneqq \bigcup_{n\leq L} \ptsLvlTP{n}_{\Psi}$,
  that is used for the evaluation of the
  non-tensor product coefficients.  
\end{remark}

The complete preprocessing phase is described by the Algorithm~\ref{algo:preprocessing}.

\begin{algorithm}
  \DontPrintSemicolon
  \Input{$\hr$}
  \BlankLine
  
  $\{ F^{n}\}_{n=0}^{L} = $ \classify{$\hr$} (Alg.~\ref{algo:classify})
  

  \For{$n=0,\dots,L$ such that $F^{n} \ne \emptyset$}{ 
    \For{$k = 1,\dots,d $}{
      compute 1D quadrature points $\ptsUnivar{n}{k} = \{ \ptbUnivar{n}{k} \}_{q_{k} = 1}^{\nptsLvlTP{n}}$ for $\bsplineUnivar{n}{k}$ (see Section~\ref{sec:WQR})\;
      
      \For{$\ell=0,\dots,n$}
      {
        compute $G_{k}^{\ell,n}$ from \eqref{eq:bsp_active_id_k}\;
      }
      compute $D_{k}^{n}$ from \eqref{eq:set_Dnk}\;
      \ForEach{$j_k \in D_k^n$}{
      $\{ \bar w_{k,j_k ,q_k}^{n} \}_{q_k \in A_{k,j_k}^{n,n}}$ = \weightsunivar{$\bsplineUnivar{n}{k}$, $j_{k}$, $\ptsUnivar{n}{k} $} (Alg.~\ref{algo:weightsunivar})\;
      }
      $ {\cal W}^{n}_k=\emptyset$\;
      \For{$\ell = 0,\dots,n$}{
        \ForEach{$i_{k} \in G_{k}^{\ell,n}$}{
          compute weights $\mathcal{W}^{\ell,n}_{k,i_k}$ from (\ref{eq:linear_comb_w_1D})\;
          $ \mathcal{W}^{n}_k = \mathcal{W}^{n}_k \cup \mathcal{W}^{\ell,n}_{k,i_k}$\;
        }
      }
    }
    compute the $d$-dimensional tensor-product points $\ptsLvlTP{n}$ from \eqref{eq:ptsTP}\;
    compute $\Psi^{n}$ from \eqref{eq:Psin}\;
    compute $\ptsLvlTP{n}_{\Psi}$ from \eqref{eq:QPsin}\;
  }
  \BlankLine
  
  \Output{$F^{n}, \mathcal{W}^{n}_1, \dots, \mathcal{W}^{n}_d, \ptsLvlTP{n}, \ptsLvlTP{n}_{\Psi}$, for $n=0,\dots,\L$ }  
  \caption{\texttt{preprocessing}}
  \label{algo:preprocessing}
\end{algorithm} 


%% file: HB_WQ_Assembly.tex
\newcommand{\supportIntersection}[0]{\Omega_{(\ell,\fnTid{i}),(\m,\fnTid{j})}}

The rows and the columns of the mass matrix are associated to the test and trial functions, respectively.
In order to emphasize the hierarchical level of a given basis function, 
we use row (or column) multi-index basis identifiers as in \eqref{eq:hr_ids}. 
Therefore the single entry of the mass matrix is denoted as 
$\left[M\right] _{(\ell,\fnTid{i}),(\m,\fnTid{j})}$ and is defined as:
\begin{equation*} \label{eq:pseudoalg_integral1}
  \left[M\right]_{(\ell,\fnTid{i}),(\m,\fnTid{j})} = 
  \int_{[0,1]^d} c(\bx) \fnLvl{\ell}{i}(\bx) \fnLvl{\m}{j}(\bx) \, d\bx \, ,
\end{equation*}
where the function $c \colon [0,1]^d \to \mathbb{R}$ incorporates the determinant of the 
Jacobian of the mapping between the parametric domain $[0,1]^d$ and the physical domain $\Omega$, 
and in general it does not have a tensor-product structure.
Recalling the quadrature rule definition (\ref{eq:rule}), and the fact that $\wtMultivar{\ell}{i} = 0$ for $\ptbMultivar{n} \notin \supp{\fnLvl{\ell}{i}}$, we 
can write
\begin{align} 
  \label{eq:mass_wq}
  \left[M\right]_{(\ell,\fnTid{i}),(\m,\fnTid{j})} \approx \quadrule{\ell}{i} (c\, \fnLvl{m}{j}) 
  = \sum_{\ptTid{q} \in \ptsTPTids{n}} \wtMultivar{\ell}{i} c(\ptbMultivar{n}) 
  \fnLvl{\m}{j}(\ptbMultivar{n}) \, d\bx,
\end{align}
where $n = \nu(\ell,\fnTid{i})$.
Using the sum-factorization approach, we exploit 
$w_{\fnTid{i},\ptTid{q}}^{n} = \prod_{k=1}^d w_{k,\fnGid{i}_{k},\ptGid{q}_{k}}^{n}$ 
and $\fnLvl{m}{j}(\ptbMultivar{n}) = \prod_{k=1}^{d} \fnUnivar{m}{k,j_k}(\bar{x}_{k,\ptGid{q}_k}^{n})$ 
and write \eqref{eq:mass_wq} in terms of nested sums:
\begin{equation}
  \label{eq:sumfact}
  \begin{aligned} 
    \quadrule{\ell}{i} (c \fnLvl{m}{j}) 
    & = \sum_{\ptGid{q}_1, \ldots, \ptGid{q}_d} \prod_{k=1}^d \left ( w_{k,\fnGid{i}_{k},\ptGid{q}_{k}}^{\ell} 
    \fnUnivar{m}{k,j_k}(\bar{x}_{k,\ptGid{q}_k}^{n}) \right ) \, c(\bar{x}_{1,\ptGid{q}_1}^{n}, \ldots, \bar{x}_{d,\ptGid{q}_d}^{n})
    \\ & = \sum_{\ptGid{q}_d} w_{d,\fnGid{i}_d,\ptGid{q}_d}^{\ell} \fnUnivar{m}{d,j_d}(\bar{x}_{d,\ptGid{q}_d}^{n})
    \left( \sum_{\ptGid{q}_{d-1}} \ldots \sum_{\ptGid{q}_1 } w_{1,\fnGid{i}_{1},\ptGid{q}_{1}}^{\ell} 
    \fnUnivar{m}{1,j_1}(\bar{x}_{1,\ptGid{q}_1}^{n}) c(\bar{x}_{1,\ptGid{q}_1}^{n}, \ldots, \bar{x}_{d,\ptGid{q}_d}^{n})\right) 
\end{aligned}
\end{equation}
where, in the summations above, each running index $\ptGid{q}_k$ ($k=1,\dots,d$)  belongs to the set
\begin{equation}
  \label{eq:pts_common_indices}
  Q_{k,i_{k},j_{k}}^{n,\ell,m} \coloneqq \{ \ptGid{q}_{k} \in \{1,\dots,R_{k}^{n} \} \colon 
  \bar{x}_{k,\ptGid{q}_{k}}^{n} \in \supp{\fnUnivar{\ell}{k,i_k}} \cap \supp{\fnUnivar{m}{k,j_k}} \} \; .
\end{equation}
Details are presented in the remaining part of this subsection, where, for the sake of notation simplicity, 
we will systematically omit the set
$Q_{k,i_{k},j_{k}}^{n,\ell,m}$ for the running index $\ptGid{q}_{k}$ in the summations.

In \eqref{eq:sumfact} we note that coefficient 
$c(\bx)$ must be evaluated at the points of the quadrature rule of
level $n$. Moreover, from \eqref{eq:prop1}--\eqref{eq:prop3} we know
that the sets $\bigl\{ F^{n} \bigr\}_{n=0}^{L}$ (excluding the empty
sets) form a partition of the hierarchical basis $\hr$. This suggests
to construct the matrix starting from an outer loop over $\bigl\{ F^{n} \bigr\}_{n=0}^{L}$, 
i.e., over the different levels of quadrature rules, then for a given level $n$, 
compute the determinant of the Jacobian at $\ptsLvlTP{n}_{\Psi}$ 
(i.e., on the points that have non-empty intersection with the support of each basis 
function in $F^n$) and set the values to be zero for the points 
$\ptsLvlTP{n} \setminus \ptsLvlTP{n}_{\Psi}$. 

The key point here is that the evaluation of the non-tensor-product coefficient $c(\bx)$ may be a costly operation, so we want to evaluate it just for the involved quadrature points, i.e., for each 
$ \ptTid{q} \in \ptsTPTids{n} $
we set:
\begin{equation}\label{eq:coeffC}
 C^{n}_{\ptTid{q}} = C^{n}_{(\ptGid{q}_{1},\dots,\ptGid{q}_{d})} \coloneqq \left\{ 
    \begin{array}{ll} 
      c(\ptbMultivar{n}) & \text{if} \; \ptbMultivar{n} \in \ptsLvlTP{n}_{\Psi} \\
      0 & \text{otherwise} 
    \end{array} \right..
\end{equation}
\begin{remark}
The coefficients $C^{n}_{\ptTid{q}}$ will be used in the innermost loop of the sum-factorization algorithm, 
so they should be stored in an efficient data structure for the data retrieval w.r.t.\ the loop index ordering 
used in the sum-factorization. 
\end{remark}

%
%

Given a quadrature level $n \in \{0,\dots,L\}$ such that $F^n \ne \emptyset$, 
we loop over $\ell,m \in \{0,\dots,n\}$ and compute the connectivity between 
the test functions of $\hrLvl{\ell} \cap F^{n}$ and the trial functions of $\hrLvl{m}$, i.e.,
\begin{equation} 
  \label{eq:connectivity}
  K_{\ell,m}^{n} \coloneqq \{ (\fnTid{i},\fnTid{j}) \in \bsplineFnTids{\ell} 
  \times \bsplineFnTids{m} \; : \; \fnLvl{\ell}{i} \in \hrLvl{\ell} \cap F^{n} 
  \, , \, \fnLvl{m}{j} \in \hrLvl{m} \, , \, \supp{\fnLvl{\ell}{i}} \cap 
  \supp{\fnLvl{m}{j}} \ne \emptyset \} \; .
\end{equation}
At this point we can apply the sum-factorization algorithm that allows
us to evaluate the mass-matrix entries.

The sum-factorization algorithm in essence is a clever way to perform the nested 
sum~\eqref{eq:sumfact}, that sequentially performs the integration along 
the directions $k=1,\ldots, d$, considering for each $k$ all pairs of indices 
$(i_k,j_k)$ that identify the weight and trial function respectively.

The integration along direction $k=1$ writes as
\begin{align}
  I^{(1)}_{(i_1),(j_1);(\ptGid{q}_2,\dots,\ptGid{q}_d)} & \coloneqq \sum_{\ptGid{q}_1} 
  w^{\ell}_{1,i_1,\ptGid{q}_1} \fnUnivar{m}{1,j_{1}}(\ptbUnivar{n}{1}) C^{n}_{(\ptGid{q}_1,\dots,\ptGid{q}_d)} \nonumber
  \\ &= \sum_{\ptGid{q}_1} w^{\ell}_{1,i_1,\ptGid{q}_1} \fnUnivar{m}{1,j_{1}}(\ptbUnivar{n}{1}) 
  I^{(0)}_{(),();(\ptGid{q}_1,\dots,\ptGid{q}_d)} \; ,   \label{eq:sf_1}
\end{align}
where we have defined $I^{(0)}_{(),();(\ptGid{q}_1,\dots,\ptGid{q}_d)} \equiv C^{n}_{(\ptGid{q}_1,\dots,\ptGid{q}_d)}$, 
which only depends on the $d$-tuple of indices associated to the quadrature points. 
Performing the summation over $\ptGid{q}_1$ we have as result $I^{(1)}_{(i_1),(j_1);(\ptGid{q}_2,\dots,\ptGid{q}_d)}$ 
that depends on the pair $(i_1,j_1)$ (related to the univariate test and trial basis along 
direction $1$) and on the $(d-1)$-tuple $(\ptGid{q}_2,\dots,\ptGid{q}_d)$ (related to univariate quadrature points 
along the directions $2,\dots,d$).
The integration along directions $k=2,\dots,d-1$ then writes as:
\begin{equation} 
  \label{eq:sf_k}
  I^{(k)}_{(i_1,\dots,i_{k}),(j_1,\dots,j_{k});(\ptGid{q}_{k+1},\dots,\ptGid{q}_d)} \coloneqq
 \sum_{\ptGid{q}_k} w^{\ell}_{k,i_k,\ptGid{q}_k} \fnUnivar{m}{k,j_{k}}(\ptbUnivar{n}{k}) I^{(k-1)}_{(i_1,\dots,i_{k-1}),(j_1,\dots,j_{k-1});(\ptGid{q}_k,\dots,\ptGid{q}_d)} \; ,
\end{equation}
and finally for $k=d$:
\begin{equation}\label{eq:sf_d}
  I^{(d)}_{(i_1,\dots,i_{d}),(j_1,\dots,j_{d});()} \coloneqq \sum_{\ptGid{q}_d} w^{\ell}_{d,i_d,\ptGid{q}_d} \fnUnivar{m}{d,j_{d}}(\ptbUnivar{n}{d}) I^{(d-1)}_{(i_1,\dots,i_{d-1}),(j_1,\dots,j_{d-1});(\ptGid{q}_d)} =  \quadrule{\ell}{i} (c \fnLvl{m}{j}) \; ,
\end{equation}
where the final expression in \eqref{eq:sf_d} is now independent of the quadrature point index but it depends on the pair of test and trial $d$-tuple $((i_1,\dots,i_{d}),(j_1,\dots,j_{d}))$ and is equal to $\quadrule{\ell}{i} (c \fnLvl{m}{j}) $.

The key point (which allows to save computations) is that the value of 
$I^{(k-1)}_{(i_1,\dots,i_{k-1}),(j_1,\dots,j_{k-1}); (\ptGid{q}_{k},\dots,\ptGid{q}_d)}$ in \eqref{eq:sf_k} 
may be needed to compute multiple values of $I^{(k)}_{(i_1,\dots,i_{k}),(j_1,\dots,j_{k});(\ptGid{q}_{k+1},\dots,\ptGid{q}_d)}$. 
To exploit this fact, when we are integrating along a direction $k$ we 
must consider all the pairs of $k$-tuples $((i_1,\dots,i_k),(j_1,\dots,j_k))$.

Accordingly, for each $k \in \{ 1,\dots, d\}$ we define the ``projection'' of the connectivity $K_{\ell,m}^n$ along the first $k$ directions:
\begin{equation} 
  \label{eq:conn_proj}
  \Pi^{(k)} K_{\ell,m}^{n} \coloneqq \bigl\{ \bigl( (i_{1},\dots,i_{k}), (j_{1},\dots,j_{k}) \bigr) 
  = \bigl( \tau_{k} \fnTid{i}, \tau_{k} \fnTid{j} \bigr) \, , \, \forall \bigl( \fnTid{i},\fnTid{j} 
  \bigr) \in K_{\ell,m}^{n} \bigr\}
\end{equation}
and then the pairs of $k$-tuple that must be considered for the efficient computation of \eqref{eq:sf_k} are just the elements of $\Pi^{(k)} K_{\ell,m}^{n}$.

The sum-factorization algorithm is then summarized by Algorithm~\ref{algo:sf}.

%
%
%

\begin{algorithm}
  \DontPrintSemicolon
  
  \Input{$\ptsTPTids{n}, \bigl\{ C^{n}_{\ptTid{q}} \bigr\}_{\ptTid{q} \in \ptsTPTids{n}}, \ell,K^{n}_{\ell,m}, \bigl\{ \mathcal{W}^{n}_k \bigr\}_{k=1}^d$}
  \BlankLine
 
  \ForEach{$(\ptGid{q}_{1},\dots,\ptGid{q}_{d}) \in \ptsTPTids{n} $}
  {
    $I^{(0)}_{(),();(\ptGid{q}_{1},\dots,\ptGid{q}_{d})} = C^{n}_{(\ptGid{q}_{1},\dots,\ptGid{q}_{d})}$\;
  }
  \For{$k=1,\dots,d$}
  {
    compute $\Pi^{(k)} K_{\ell,m}^{n}$ from (\ref{eq:conn_proj})\;
	\ForEach{$(i_{k},j_{k}) \in \{ \{1,\dots,\ndofsUnivar{\ell}{k}\} \times \{1,\dots,\ndofsUnivar{m}{j}\} 
	\, : \, \supp{\fnUnivar{\ell}{k,i_k}} \cap \supp{\fnUnivar{m}{k,j_k}} \ne \emptyset \}$}
	{
	  compute $Q_{k,i_{k},j_{k}}^{n,\ell,m}$ from \eqref{eq:pts_common_indices};
	}
    \ForEach{$\bigl( (i_{1},\dots,i_{k}), (j_{1},\dots,j_{k}) \bigr) \in \Pi^{(k)} K_{\ell,m}^{n}$}
    {
      retrieve $\mathcal{W}^{\ell,n}_{k,i_k}$ from $\mathcal{W}^{n}_{k}$ (see~(\ref{eq:linear_comb_w_1D}))\;
      \ForEach{$(\ptGid{q}_{k+1},\dots,\ptGid{q}_{d}) \in \{1,\dots,\nptsUnivar{n}{k+1} \} \times \dots \times \{1,\dots,\nptsUnivar{n}{d} \} $}
      {
        compute $I^{(k)}_{(i_{1},\dots,i_{k}),(j_{1},\dots,j_{k});(\ptGid{q}_{k+1},\dots,\ptGid{q}_{d})}$ from \eqref{eq:sf_k}\;
      }
    }

  }

  \BlankLine
  \Output{$I^{(d)} \equiv \bigl\{ I^{(d)}_{(\fnTid{i}),(\fnTid{j});()} \bigr\}_{(\fnTid{i},\fnTid{j}) \in K_{\ell,m}^n}$}
  \caption{\texttt{sum\_factorization}}
  \label{algo:sf}
\end{algorithm}

The algorithm for the matrix formation  is depicted by
Algorithm~\ref{algo:assemble}.

\begin{algorithm}
  \DontPrintSemicolon
  {
    \Input{$\hr, \{ F^{n}, \ptsLvlTP{n}, \ptsLvlTP{n}_{\Psi}, \mathcal{W}^{n}_1, \dots, \mathcal{W}^{n}_d \}_{n=0}^{L}$,  $c$}

    \BlankLine
    
    \ForEach{$n \in \{0,\dots,L \}$ such that $F^{n} \ne
      \emptyset$} 
    {
      compute $\bigl\{ C^{n}_{\ptTid{q}} \bigr\}_{\ptTid{q} \in \ptsTPTids{n}}$ from (\ref{eq:coeffC})\;
      
	  \ForEach{$\ell \in \{ 0,\dots,n \}$ such that  $\hrLvl{\ell} \cap F^{n} \ne \emptyset$}
	  {
	    \ForEach{$m \in \{ 0,\dots,n \}$}
	    {
          compute $K_{\ell,m}^{n}$ from (\ref{eq:connectivity})\;
          
%
%
%
%

          
            $ \bigl\{ [M]_{(\ell,\fnTid{i}),(m,\fnTid{j})} \bigr\}_{(\fnTid{i},\fnTid{j}) \in K_{\ell,m}^n} =$ \texttt{sum\_factorization($\ptsTPTids{n},\bigl\{ C^{n}_{\ptTid{q}} \bigr\}_{\ptTid{q} \in \ptsTPTids{n}},\ell,K^{n}_{\ell,m},\bigl\{ \mathcal{W}^{n}_k \bigr\}_{k=1}^d$)} (Alg.~\ref{algo:sf})\;

        }
      }
    }
  }
  
  \BlankLine
  \Output{$M$}
  \caption{\texttt{compute\_matrix}}
  \label{algo:assemble}
\end{algorithm}


%% file: Cost.tex
We now want to estimate the total computational cost of the matrix  formation. There are mainly three steps that contribute to this cost: the evaluation of the non-tensor product coefficient $c$, the computation of the weights, and the computation of the matrix entries via sum-factorization. 

The coefficient $c$ has to be evaluated for every active quadrature point. Quadrature points are more dense for elements that are adjacent to the boundary of $\Omega$. However the total number of active quadrature points is dominated from the interior part in all cases of interest. Recalling \eqref{eq:hier_pointSet} and $\nu-\ell \leq r-1$,  the number of quadrature points that belong to interior  elements is bounded by
 
\begin{align} \label{eq:quad_points_bound_NH} 
\sum_{ (\ell,\fnTid{i}) \in \mathcal{I}_{\mathcal{H}}} \#  \ptsLvlTP{\left(\ell,\fnTid{i}\right)}  \leq \sum_{ (\ell,\fnTid{i}) \in \mathcal{I}_{\mathcal{H}}}   (2^{r-1}(p+1))^{d} = O( 2^{dr}p^d  \ndofsHr ).
  \end{align}
We remark that  bound above is not sharp especially  for what concerns its dependence on $p$, since quadrature points in different $\ptsLvlTP{\left(\ell,\fnTid{i}\right)}$ may coincide.



As for the computation of the weights, we recall that we have to solve a system of the form \eqref{eq:exactness1D} for every univariate index $j_k$, for $k=1,\ldots,d$. Since the number of univariate indices is bounded by the number of multi-indices $\ndofsHr$, and since each of these linear system has $O(p)$ unknown nonzero weights 
the cost to compute them all using a direct solver is bounded by $ O(p^3 \ndofsHr) $ flops. 

If we compare the bound on this cost with the one related to the computation of the matrix entries (derived below), we see that they have the same order with respect to $p$ for $d=2$ and that the former has lower order for $d=3$. Note also that this bound does not depend on the admissibility parameter $r$.


%
%

We now discuss the computation of the matrix entries. Following the structure of Algorithm~\ref{algo:assemble}, we fix $n \in \left\lbrace 0, \ldots , L \right\rbrace  $ and $ \ell, m \in \left\lbrace 0, \ldots, n \right\rbrace $ and consider the computation of the matrix entries \eqref{eq:mass_wq} for all $(\fnTid{i},\fnTid{j}) \in K_{\ell,m}^{n}$, as performed by Algorithm \ref{algo:sf}. 

%

As a preliminary step, we observe that for any fixed direction $k \in \lbrace 1,\ldots,d \rbrace$ and any fixed index value $i_k$, the number of indices  $j_k$ that must be considered in \eqref{eq:sumfact} is clearly bounded by the number of basis functions of level $m$ whose support intersects the support of $\fnUnivar{\ell}{k,i_k}$. It can be verified that the latter number is bounded by $2p + 1$ when $m \leq \ell$, and by $ 2^{m - \ell} (p+1)+p $ when $m > \ell$. In both cases, this number is bounded by $2^{n - \ell +1} (p + 1)$, since $n \geq \max \lbrace m, \ell \rbrace$.

Moreover, 
again for any fixed direction $k$ and index value $i_k$,
the active quadrature points $\ptbUnivar{n}{k}$ are the ones belonging to the support of $\fnUnivar{\ell}{k,i_k}$; since we have $2$ quadrature points on each interior element of level $n$, or $p+1$ on the elements that touch the boundary, and the support of $\fnUnivar{\ell}{k,i_k}$ contains at most $2^{n - \ell } (p + 1)$ elements of level $n$, we conclude that there are at most $2^{n - \ell + 1} (p + 1)$ active quadrature points if  $\fnUnivar{\ell}{k,i_k}$  does not touch the boundary, or at most   $\left( 2^{n - \ell + 1} + 1 \right) (p + 1) $ quadrature points if $\fnUnivar{\ell}{k,i_k}$   touches the boundary. Typically, the cost is dominated by the quadrature at the interior, therefore we assume that the number of index values taken by $\ptGid{q}_k$ in the $k$-th sum of  \eqref{eq:sumfact} is roughly $2^{n - \ell + 1} (p + 1)$. 

We are now ready to estimate the cost of computing \eqref{eq:sumfact}. As a first step, we evaluate the innermost sum \eqref{eq:sf_1} for all relevant values of $i_1$, $j_1$ and $\ptGid{q}_2, \ldots,\ptGid{q}_d$. Of course in the sum we only need to consider the nonzero terms, and we observe that the term corresponding to a fixed $\ptGid{q}_1$ is nonzero only for the $p+1$ values of the index $j_1$ such that $\fnUnivar{m}{1,j_1}(\bar{x}_{\ptGid{q}_1}^{n}) \neq 0$. 
%
%
{
Note that if we preliminary multiply  $w_{1,\fnGid{i}_{1},\ptGid{q}_{1}}^{\ell} \fnUnivar{m}{1,j_1}(\bar{x}_{\ptGid{q}_1}^{n})$ for all such values of $\ptGid{q}_1$ and $j_1$ (which has a negligible cost), the compution of the sum \eqref{eq:sf_1} requires 2 flops for each of its nonzero terms.

Since each index $\ptGid{q}_1,\ldots,\ptGid{q}_d$, can take up to $2^{n - \ell + 1} (p + 1)$ values, and the number of values taken by $i_1 = \tau_{1} \fnTid{i}$ is bounded by the number of multi-indices $\fnTid{i}$ belonging to $F^{n} \cap \hrLvl{\ell}$, }
the cost of the first step is bounded by
\begin{align} \label{eq:cost_inner_sum} 2 \left( p+1 \right)^{d+1} \, 2^{d(n - \ell + 1)}  N_{n,\ell} \text{ flops}. \end{align}
{
where
$$ N_{n,\ell} := \vert F^{n} \cap \hrLvl{\ell}\vert. $$


For $k = 2,\ldots,d-1$, the $k$-th step of the sum-factorization requires the computation of \eqref{eq:sf_k}
%
%
for all values of $i_1, \ldots, i_{k}$, $j_1, \ldots, j_{k}$ and $\ptGid{q}_{k+1}, \ldots, \ptGid{q}_d$, where the inner sum $I^{(k-1)}_{(i_1,\dots,i_{k-1}),(j_1,\dots,j_{k-1});(\ptGid{q}_k,\dots,\ptGid{q}_d)}$ has already been computed for all the relevant index values. 

Since $\left( i_1, \ldots, i_{k}\right) = \tau_{k} \fnTid{i}$ the total number of 
$k$-tuples $\left( i_1, \ldots, i_{k}\right)$ that have to be considered is again bounded by the number of multi-indices $N_{n,\ell}$.
}

Moreover, again we observe that for each value of $\ptGid{q}_k$ there are only $p+1$ values of $j_k$ that contribute to the sum, and since the number of values taken by each index $\ptGid{q}_{k}, \ldots, \ptGid{q}_d$ and $j_1, \ldots, j_{k-1}$ is bounded by $2^{n - \ell + 1} (p + 1)$, the cost of this step is again bounded by \eqref{eq:cost_inner_sum}. 
{ With similar arguments, it can be shown that this is true also for the $d-$th step of the sum-factorization \eqref{eq:sf_k}.}

We conclude that the cost of the whole sum-factorization step is bounded by
\begin{align*} 2d \left( p+1 \right)^{d+1} \, 2^{d(n - \ell + 1)}  N_{n,\ell} \text{ flops}. \end{align*}

We sum the above expression for all values of $n, \ell$ and $m$, and observe that for a fixed level $\ell$ the number of levels $m$ that interact with it is at most $2r-1$. Thus, a bound on the total cost for the matrix entries computation is given by
\begin{align}\label{eq:total_cost} 2d (2r-1) \left( p+1 \right)^d \sum_{n} \sum_{\ell \leq n}  2^{d(n - \ell + 1)} N_{n,\ell} \text{ flops}. \end{align}

We can derive a more explicit bound on the cost of the matrix entries computation if take a further step and observe that $n - \ell + 1 \leq r$ and that
$$ \sum_{n} \sum_{\ell \leq n} N_{n,\ell} = \ndofsHr $$
Hence the total cost for the matrix entries computation \eqref{eq:total_cost} is bounded by
\begin{align}\label{eq:total_cost_bound} 2d (2r-1) 2^{dr} \left( p+1 \right)^{d+1} \ndofsHr = O \left( dr 2^{dr} p^{d+1} \ndofsHr \right) \text{ flops}. \end{align}

We observe that, similarly as in the bound on the active quadrature points \eqref{eq:quad_points_bound_NH}, the latter expression grows exponentially with respct to the admissibilty parameter $r$, and this effect worsen with the increasing of the dimension $d$. This might seem unsatisfactory, but we emphasize that \eqref{eq:total_cost_bound} is easily a rather pessimistic bound. Indeed, a careful analysis of the derivation of \eqref{eq:total_cost_bound} reveals that we are essentially assuming that every hierarchical B-spline basis function $\fnLvl{\ell}{i} $, with  $(\ell,\fnTid{i}) \in \dofsHr$, interacts with all the admissible levels. 
In many practical cases, however, refinement is perfomred only in specific regions of the domain, e.g., in the neighbourhood of  low dimensional manifolds, and as a result the number of basis functions that interact with all the admissible levels is limited.



%% file: HB_WQ_NumericalTests.tex
{

The numerical tests comprise of the $L^{2}$-projection of the function $f \colon \Omega \to \mathbb{R}$,
\begin{equation} \label{eq:func_to_project}
  f(\mathbf{x}) = \exp \biggl[ - \Bigl( \frac{\Vert \mathbf{x} - \mathbf{x}_{0} \Vert - 1}{\beta} \Bigr)^{2} \biggr],
\end{equation}
where the physical domain $\Omega \in \RR^{d}$, the parameter $\beta \in \RR$ and the point $\mathbf{x}_{0} \in \mathbb{R}^{d}$ are specified below for the cases $d=2$ and $d=3$, using a $\admissibility$-admissible hierarchical B-spline basis (see Section~\ref{sec:HB}) of degree $p$, with different values for the admissibility parameter $\admissibility$.

For each value of the admissibility parameter $\admissibility$, a nested sequence of hierarchical B-spline spaces is constructed \cite{bracco2018b}. The adaptive mesh refinement is steered by the ``error estimator'', which is simply the $L^2$-error between the computed $L^{2}$ projection and the function \eqref{eq:func_to_project} and by using the D\"orfler marking strategy \cite{dorfler1996} with parameter $\theta_{\ast}=0.2$.

For each refinement step, we perform a simulation using the standard element-base Gaussian quadrature (using $p+1$ quadrature point along each direction of the element) to build the mass matrix (and right hand side) and then, using the same sequence of hierarchical spaces we compute the mass-matrix using the proposed hierarchical WQ algorithm.

\begin{remark}
All the numerical tests were performed using the IGATOOLS library \cite{igatools2015}, on a single core of an Intel Xeon E5-2470 processor running at 2.3~GHz. 
In order to alleviate the random fluctuations in the elapsed CPU time, 
all plots involving CPU time refer to the average CPU time of multiple (5 for the 2D case and 3 for the 3D case) runs of the same simulation.
\end{remark}

\begin{remark}
  In all plots the lower limit of the CPU time is set to $10^{-1}$ seconds to reduce the effect of random time fluctuations, due to the CPU scheduling. 
\end{remark}

\subsection{2D case}
For this case the physical domain $\Omega \subset \RR^{2}$ is defined as the image of the two-dimensional parametric domain 
$\hat{\Omega} = [1,2] \times [\frac{\pi}{4},\frac{3\pi}{4}]$
through the (polar) map
\begin{equation*}
  F(\rho,\theta) = 
  \begin{pmatrix}
    \rho \cos{\theta}  \\ 
    \rho \sin{\theta}
  \end{pmatrix} \, ,
\end{equation*}
while $\beta = 5 \cdot 10^{-3}$ and $\mathbf{x}_{0} = (0,\tfrac{5}{2})$ (see Figure~\ref{fig:2Dp2adm2ref51}).

\begin{figure}  
  \centerline{\includegraphics[width=\textwidth]{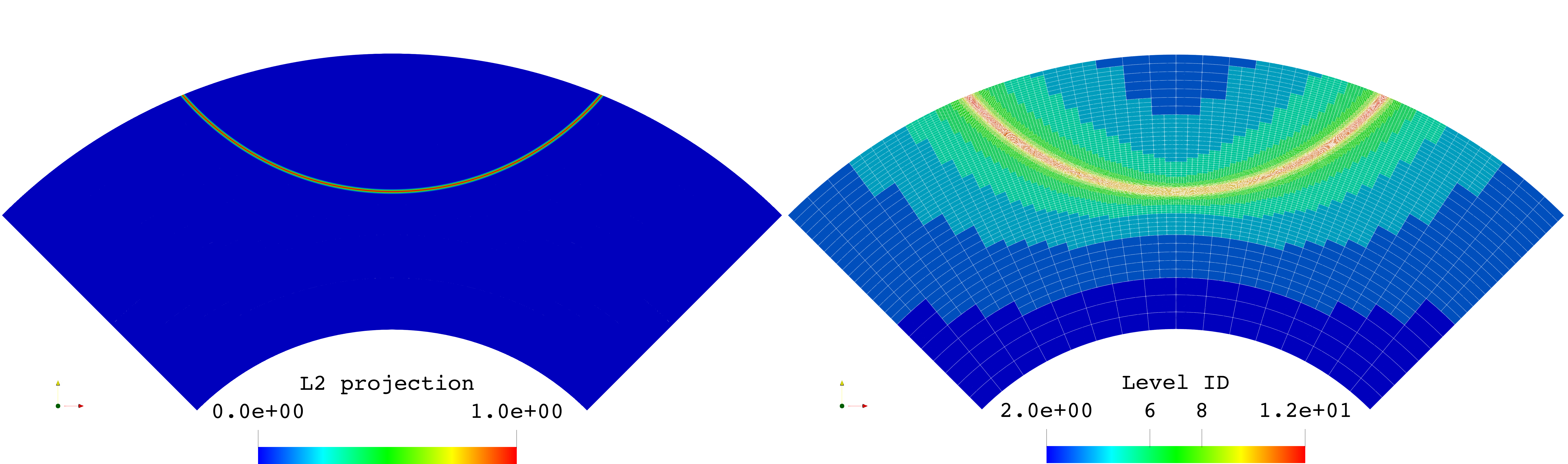}}
  \caption{$L^{2}$-projection (left) and element levels (right) after 51 adaptive refinements 
  for the 2D case with degree $p=2$ and admissibility $\admissibility=2$. For this configuration 
  the space contains 867947 degrees of freedom, and the $L^{2}$-error between the function (\refeq{eq:func_to_project}) 
  and its $L^{2}$-projection is $\approx \num{3.4e-7}$.
  \label{fig:2Dp2adm2ref51}}
\end{figure}

For this case we performed simulations using the admissibility parameters $\admissibility = 2,3$ 
and for each value of $\admissibility$ we used the degrees $p=2,\dots,6$.

Regarding the cases with $r=2$ we can observe from the plots in Figure~\ref{fig:dofs_vs_total_assembly_time_2D_adm=2} that the total time (preprocessing $+$ matrix computation) w.r.t.\ $\ndofsHr$   or the WQ approach seems to be nearly independent from the degree $p$, while for the element-based Gaussian approach we note that the cost increases with $p$ (as expected). 
Moreover, also the most favorable case for the element-based Gaussian approach (i.e., $p=2$) costs more of any of the WQ cases we have tested (except for some specific space configurations when $p=3$). As a result, we can conclude that if one wants achieve a very low error level ($< 10^{-7}$), the best strategy in terms of CPU time needed to build the matrix is to use WQ with high degree (see Figure~\ref{fig:err_vs_total_assembly_time_2D_adm=2}).

Regarding the CPU cost of the WQ approach, in Figure~\ref{fig:err_vs_time_WQ_details_2D_adm=2} are shown (for the degrees $p=2,\dots,5$) the preprocessing cost (Algorithm~\ref{algo:preprocessing}) and the matrix computation cost (Algorithm~\ref{algo:assemble}), that is split in the time needed to evaluate the coefficients in \eqref{eq:coeffC} (for $n=0,\dots,L$) and the rest of the algorithm (i.e., the computation of the connectivities $K_{\ell,m}^{n}$ from \eqref{eq:connectivity} and the sum-factorization). 
From the plots in Figure~\ref{fig:err_vs_time_WQ_details_2D_adm=2} we have that the asymptotic behaviour of the costs is the same for all different degrees, resulting in the dominant cost being the formation of the matrix whereas the cost for the preprocessing is smaller but not negligible (at least for the tested cases). It is worthy to note that for low number of degrees of freedom, the main cost is due to the preprocessing. Moreover, the cost for evaluating the coefficients in in Eq.~\eqref{eq:coeffC} depends on the number of points in $\{ \ptsLvlTP{n}_{\Psi} \}_{n=0}^{L}$ and che cost of evaluation of the function $c$ at a single point. In our case $c$ is just the determinant of the Jacobian of the mapping, resulting in low CPU time w.r.t.\ the other two main costs.

Regarding the case with admissibility parameter $\admissibility = 3$, we observe from the plots in Figure~\ref{fig:dofs_vs_total_assembly_time_2D_adm=3} that both approaches (element-based Gaussian quadrature and WQ) have an higher cost  (for a given number of dofs) for all tested degrees w.r.t. the case with $\admissibility = 2$ (Figure~\ref{fig:dofs_vs_total_assembly_time_2D_adm=2}), but the WQ approach seems to be nearly independent of the degree and it less expensive w.r.t. the element-based Gaussian approach of degree $>=3$.

\input{numerical_tests_final/HB_WQ_tests_2D_adm2_dofs_vs_total_assembly_time.tex}
\input{numerical_tests_final/HB_WQ_tests_2D_adm2_L2error_vs_total_assembly_time.tex}
\input{numerical_tests_final/HB_WQ_tests_2D_adm2_dofs_vs_WQ_details.tex}

\input{numerical_tests_final/HB_WQ_tests_2D_adm3_dofs_vs_total_assembly_time.tex}

\subsection{3D case}
For this case the physical domain $\Omega \subset \RR^{3}$ is defined as the image of the two-dimensional parametric domain 
$\hat{\Omega} = [1,2] \times [\frac{\pi}{4},\frac{3\pi}{4}] \times [0,\frac{\pi}{2}]$
through the (polar) map
\begin{equation*}
  F(\rho,\theta,\phi) = 
  \begin{pmatrix}
    \rho \cos{\theta}  \\ 
    \rho \sin{\theta} \cos{\phi} \\
    \rho \sin{\theta} \sin{\phi}
  \end{pmatrix} \, ,
\end{equation*}
while $\beta = 0.1$ and $\mathbf{x}_{0} = (0,\tfrac{5}{2},0)$ (see the Figure~\ref{fig:3Dp2adm2ref27}).

\begin{figure} 
  \centerline{\includegraphics[width=\textwidth]{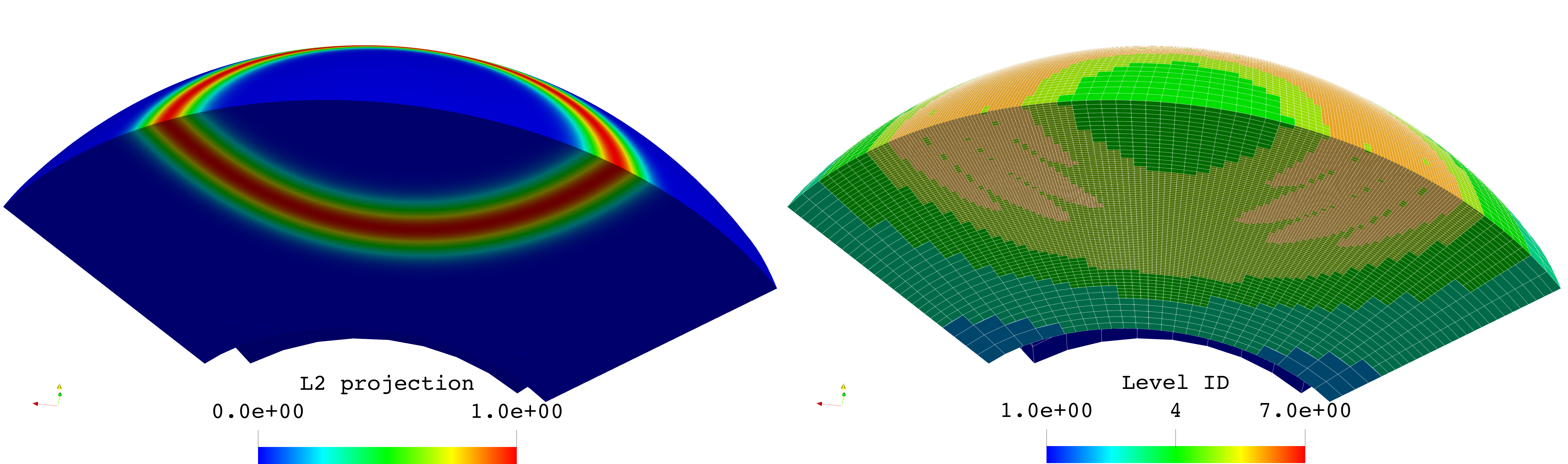}}
  \caption{$L^{2}$-projection (left) and element levels (right) after 27 adaptive refinements 
  for the 3D case with degree $p=2$ and admissibility $\admissibility=2$. For this configuration 
  the space contains 754614 dofs, and the $L^{2}$-error between the function (\refeq{eq:func_to_project}) 
  and its $L^{2}$-projection is $\approx \num{8.5e-6}$.
  \label{fig:3Dp2adm2ref27}.}
\end{figure}

For this case we performed simulations using the admissibility parameter $\admissibility = 2$ 
and the degrees $p=2,\dots,5$.

We can observe from the plots in Figures~\ref{fig:dofs_vs_total_assembly_time_3D_adm=2} and \ref{fig:err_vs_total_assembly_time_3D_adm=2} that the WQ approach outperforms the element-based Gaussian approach. In fact, considering the CPU time
w.r.t.\ $\ndofsHr$, the total time (preprocessing $+$ matrix computation) for the WQ approach seems to be mildly dependent from the degree $p$, while for the element-based Gaussian approach we note that the cost increases with $p$, by a factor higher than the 2D case (as expected). 
Moreover, also the most favorable case for the element-based Gaussian approach (i.e. $p=2$) costs more of any of the WQ cases we have tested. As result, we have that if one want achieve a low error level, the best strategy in terms of CPU time needed to build the matrix is to use WQ with high degree (see Figure~\ref{fig:err_vs_total_assembly_time_3D_adm=2}).

Regarding the CPU cost of the WQ approach, in Figure~\ref{fig:err_vs_time_WQ_details_3D_adm=2} are shown (for the degrees $p=2,\dots,5$) the preprocessing cost (Algorithm~\ref{algo:preprocessing}) and the matrix computation cost (Algorithm~\ref{algo:assemble}), that is split in the time needed to evaluate the coefficients in Eq.~\eqref{eq:coeffC} (for $n=0,\dots,L$) and the rest of the algorithm (i.e. the computation of the connectivities $K_{\ell,m}^{n}$ from \eqref{eq:connectivity} and the sum-factorization). 
In this case, w.r.t. the 2D case we observe that for the degree $p=2$, when we have a number of degrees of freedom $< \num{1.5e5}$, the dominant cost can be attributed to the evaluation of the coefficients in \eqref{eq:coeffC} whereas for an higher number of degrees of freedom, the dominant cost is due to the sum-factorization$+$computation of the connectivities $K_{\ell,m}^{n}$ phase. When we increase the degree we observe that this cost becomes the dominant one when $p>3$ (and conversely the cost for the coefficients evaluation becomes the lowest of the three costs for $p>3$).

\input{numerical_tests_final/HB_WQ_tests_3D_adm2_dofs_vs_total_assembly_time.tex}
\input{numerical_tests_final/HB_WQ_tests_3D_adm2_L2error_vs_total_assembly_time.tex}
\input{numerical_tests_final/HB_WQ_tests_3D_adm2_dofs_vs_WQ_details.tex}

} 

%% file: numerical_tests_final/HB_WQ_tests_2D_adm2_dofs_vs_total_assembly_time.tex

\pgfplotstableread[skip first n=1]
{./numerical_tests_final/time_and_errors_projection_l2_2D_adm_2_degree_2_func_0_beta=0.005_threads_1_adaptive_refinement_wqMat_gaussRHS.dat}\wqdataDegTwo

\pgfplotstableread[skip first n=1]{./numerical_tests_final/time_and_errors_projection_l2_2D_adm_2_degree_2_func_0_beta=0.005_threads_1_adaptive_refinement_gaussMat_gaussRHS.dat}\gaussdataDegTwo

\pgfplotstableread[skip first n=1]
{./numerical_tests_final/time_and_errors_projection_l2_2D_adm_2_degree_3_func_0_beta=0.005_threads_1_adaptive_refinement_wqMat_gaussRHS.dat}\wqdataDegThree

\pgfplotstableread[skip first n=1]{./numerical_tests_final/time_and_errors_projection_l2_2D_adm_2_degree_3_func_0_beta=0.005_threads_1_adaptive_refinement_gaussMat_gaussRHS.dat}\gaussdataDegThree

\pgfplotstableread[skip first n=1]
{./numerical_tests_final/time_and_errors_projection_l2_2D_adm_2_degree_4_func_0_beta=0.005_threads_1_adaptive_refinement_wqMat_gaussRHS.dat}\wqdataDegFour

\pgfplotstableread[skip first n=1]{./numerical_tests_final/time_and_errors_projection_l2_2D_adm_2_degree_4_func_0_beta=0.005_threads_1_adaptive_refinement_gaussMat_gaussRHS.dat}\gaussdataDegFour

\pgfplotstableread[skip first n=1]
{./numerical_tests_final/time_and_errors_projection_l2_2D_adm_2_degree_5_func_0_beta=0.005_threads_1_adaptive_refinement_wqMat_gaussRHS.dat}\wqdataDegFive

\pgfplotstableread[skip first n=1]{./numerical_tests_final/time_and_errors_projection_l2_2D_adm_2_degree_5_func_0_beta=0.005_threads_1_adaptive_refinement_gaussMat_gaussRHS.dat}\gaussdataDegFive

\pgfplotstableread[skip first n=1]
{./numerical_tests_final/time_and_errors_projection_l2_2D_adm_2_degree_6_func_0_beta=0.005_threads_1_adaptive_refinement_wqMat_gaussRHS.dat}\wqdataDegSix

\pgfplotstableread[skip first n=1]{./numerical_tests_final/time_and_errors_projection_l2_2D_adm_2_degree_6_func_0_beta=0.005_threads_1_adaptive_refinement_gaussMat_gaussRHS.dat}\gaussdataDegSix

\begin{center}

\tikzset{every mark/.append style={scale=0.5}}

\begin{figure}[!h]
\begin{tikzpicture}
\begin{loglogaxis}[
    width=\textwidth,
    height=0.7\textwidth,
    title={2D case: DOFs vs. Total time for the matrix computation, admissibility $\admissibility=2$.},
    xlabel={$\ndofsHr$},
    ylabel={Time (sec.)},
    xmin=400,xmax=2.5e6,
    ymin=1e-1,ymax=1e3,
    legend pos=north west,
    legend columns=2, 
    legend style={
      /tikz/column 2/.style={column sep=15pt,},
    },
    grid=major
  ]
    


  \addplot[red,mark=*,mark size=1pt]
    table [
      x index=2,
      y expr=\thisrowno{11}+\thisrowno{15}
      ] \wqdataDegTwo;
   \addlegendentry{WQ-HR $p=2$}

  \addplot[red,mark=*,style=densely dotted]
    table [
      x index=2,
      y index=15
      ] \gaussdataDegTwo;
   \addlegendentry{Gauss $p=2$}



  \addplot[blue,mark=*,mark size=1pt]
    table [
      x index=2,
      y expr=\thisrowno{11}+\thisrowno{15}
      ] \wqdataDegThree;
   \addlegendentry{WQ-HR $p=3$}

  \addplot[blue,mark=*,style=densely dotted]
    table [
      x index=2,
      y index=15
      ] \gaussdataDegThree;
   \addlegendentry{Gauss $p=3$}



  \addplot[magenta,mark=*,mark size=1pt]
    table [
      x index=2,
      y expr=\thisrowno{11}+\thisrowno{15}
      ] \wqdataDegFour;
   \addlegendentry{WQ-HR $p=4$}

  \addplot[magenta,mark=*,style=densely dotted]
    table [
      x index=2,
      y index=15
      ] \gaussdataDegFour;
   \addlegendentry{Gauss $p=4$}




  \addplot[black,mark=*,mark size=1pt]
    table [
      x index=2,
      y expr=\thisrowno{11}+\thisrowno{15}
      ] \wqdataDegFive;
   \addlegendentry{WQ-HR $p=5$}

  \addplot[black,mark=*,style=densely dotted]
    table [
      x index=2,
      y index=15
      ] \gaussdataDegFive;
   \addlegendentry{Gauss $p=5$}



  \addplot[cyan,mark=*,mark size=1pt]
    table [
      x index=2,
      y expr=\thisrowno{11}+\thisrowno{15}
      ] \wqdataDegSix;
   \addlegendentry{WQ-HR $p=6$}

  \addplot[cyan,mark=*,style=densely dotted]
    table [
      x index=2,
      y index=15
      ] \gaussdataDegSix;
   \addlegendentry{Gauss $p=6$}


\end{loglogaxis}
\end{tikzpicture}
\caption{$\ndofsHr$ vs.\ total time for the matrix computation for the 2D case and admissibility parameter $\admissibility=2$. For WQ-HR the total time for the matrix computation is the sum of the time needed for the preprocessing (Section~\ref{sec:preprocessing}) and for the matrix formation (Section~\ref{sec:formation_algo}).}
\label{fig:dofs_vs_total_assembly_time_2D_adm=2}
\end{figure}
%
%

\end{center}

%% file: numerical_tests_final/HB_WQ_tests_2D_adm2_L2error_vs_total_assembly_time.tex

\pgfplotstableread[skip first n=1]
{./numerical_tests_final/time_and_errors_projection_l2_2D_adm_2_degree_2_func_0_beta=0.005_threads_1_adaptive_refinement_wqMat_gaussRHS.dat}\wqdataDegTwo

\pgfplotstableread[skip first n=1]{./numerical_tests_final/time_and_errors_projection_l2_2D_adm_2_degree_2_func_0_beta=0.005_threads_1_adaptive_refinement_gaussMat_gaussRHS.dat}\gaussdataDegTwo

\pgfplotstableread[skip first n=1]
{./numerical_tests_final/time_and_errors_projection_l2_2D_adm_2_degree_3_func_0_beta=0.005_threads_1_adaptive_refinement_wqMat_gaussRHS.dat}\wqdataDegThree

\pgfplotstableread[skip first n=1]{./numerical_tests_final/time_and_errors_projection_l2_2D_adm_2_degree_3_func_0_beta=0.005_threads_1_adaptive_refinement_gaussMat_gaussRHS.dat}\gaussdataDegThree

\pgfplotstableread[skip first n=1]
{./numerical_tests_final/time_and_errors_projection_l2_2D_adm_2_degree_4_func_0_beta=0.005_threads_1_adaptive_refinement_wqMat_gaussRHS.dat}\wqdataDegFour

\pgfplotstableread[skip first n=1]{./numerical_tests_final/time_and_errors_projection_l2_2D_adm_2_degree_4_func_0_beta=0.005_threads_1_adaptive_refinement_gaussMat_gaussRHS.dat}\gaussdataDegFour

\pgfplotstableread[skip first n=1]
{./numerical_tests_final/time_and_errors_projection_l2_2D_adm_2_degree_5_func_0_beta=0.005_threads_1_adaptive_refinement_wqMat_gaussRHS.dat}\wqdataDegFive

\pgfplotstableread[skip first n=1]{./numerical_tests_final/time_and_errors_projection_l2_2D_adm_2_degree_5_func_0_beta=0.005_threads_1_adaptive_refinement_gaussMat_gaussRHS.dat}\gaussdataDegFive

\pgfplotstableread[skip first n=1]
{./numerical_tests_final/time_and_errors_projection_l2_2D_adm_2_degree_6_func_0_beta=0.005_threads_1_adaptive_refinement_wqMat_gaussRHS.dat}\wqdataDegSix

\pgfplotstableread[skip first n=1]{./numerical_tests_final/time_and_errors_projection_l2_2D_adm_2_degree_6_func_0_beta=0.005_threads_1_adaptive_refinement_gaussMat_gaussRHS.dat}\gaussdataDegSix

\begin{center}

\tikzset{every mark/.append style={scale=0.5}}

\begin{figure}[!h]
\begin{tikzpicture}
\begin{loglogaxis}[
    width=\textwidth,
    height=0.7\textwidth,
    title={2D case: $L^{2}$-error vs. Total time for the matrix computation, admissibility $\admissibility=2$.},
    xlabel={$L^{2}$-error},
    ylabel={Time (sec.)},
    xmin=2e-10,xmax=1.5e-1,
    ymin=1e-1,ymax=1e3,
    legend pos=south west,
    legend columns=2, 
    legend style={
      /tikz/column 2/.style={column sep=15pt,},
    },
    grid=major
  ]
    


  \addplot[red,mark=*,mark size=1pt]
    table [
      x index=4,
      y expr=\thisrowno{11}+\thisrowno{15}
      ] \wqdataDegTwo;
   \addlegendentry{WQ-HR $p=2$}

  \addplot[red,mark=*,style=densely dotted]
    table [
      x index=4,
      y index=15
      ] \gaussdataDegTwo;
   \addlegendentry{Gauss $p=2$}



  \addplot[blue,mark=*,mark size=1pt]
    table [
      x index=4,
      y expr=\thisrowno{11}+\thisrowno{15}
      ] \wqdataDegThree;
   \addlegendentry{WQ-HR $p=3$}

  \addplot[blue,mark=*,style=densely dotted]
    table [
      x index=4,
      y index=15
      ] \gaussdataDegThree;
   \addlegendentry{Gauss $p=3$}



  \addplot[magenta,mark=*,mark size=1pt]
    table [
      x index=4,
      y expr=\thisrowno{11}+\thisrowno{15}
      ] \wqdataDegFour;
   \addlegendentry{WQ-HR $p=4$}

  \addplot[magenta,mark=*,style=densely dotted]
    table [
      x index=4,
      y index=15
      ] \gaussdataDegFour;
   \addlegendentry{Gauss $p=4$}




  \addplot[black,mark=*,mark size=1pt]
    table [
      x index=4,
      y expr=\thisrowno{11}+\thisrowno{15}
      ] \wqdataDegFive;
   \addlegendentry{WQ-HR $p=5$}

  \addplot[black,mark=*,style=densely dotted]
    table [
      x index=4,
      y index=15
      ] \gaussdataDegFive;
   \addlegendentry{Gauss $p=5$}



  \addplot[cyan,mark=*,mark size=1pt]
    table [
      x index=4,
      y expr=\thisrowno{11}+\thisrowno{15}
      ] \wqdataDegSix;
   \addlegendentry{WQ-HR $p=6$}

  \addplot[cyan,mark=*,style=densely dotted]
    table [
      x index=4,
      y index=15
      ] \gaussdataDegSix;
   \addlegendentry{Gauss $p=6$}

\end{loglogaxis}
\end{tikzpicture}
\caption{$L^{2}$-error vs.\ total time for the matrix computation for the 2D case and admissibility parameter $\admissibility=2$. For WQ-HR the total time for the matrix computation is the sum of the time needed for the preprocessing (Section~\ref{sec:preprocessing}) and for the matrix formation (Section~\ref{sec:formation_algo}).}
\label{fig:err_vs_total_assembly_time_2D_adm=2}
\end{figure}
%
%

\end{center}

%% file: numerical_tests_final/HB_WQ_tests_2D_adm2_dofs_vs_WQ_details.tex

\pgfplotstableread[skip first n=1]
{./numerical_tests_final/time_and_errors_projection_l2_2D_adm_2_degree_2_func_0_beta=0.005_threads_1_adaptive_refinement_wqMat_gaussRHS.dat}\wqdataDegTwo

\pgfplotstableread[skip first n=1]
{./numerical_tests_final/time_and_errors_projection_l2_2D_adm_2_degree_3_func_0_beta=0.005_threads_1_adaptive_refinement_wqMat_gaussRHS.dat}\wqdataDegThree

\pgfplotstableread[skip first n=1]
{./numerical_tests_final/time_and_errors_projection_l2_2D_adm_2_degree_4_func_0_beta=0.005_threads_1_adaptive_refinement_wqMat_gaussRHS.dat}\wqdataDegFour

\pgfplotstableread[skip first n=1]
{./numerical_tests_final/time_and_errors_projection_l2_2D_adm_2_degree_5_func_0_beta=0.005_threads_1_adaptive_refinement_wqMat_gaussRHS.dat}\wqdataDegFive

\pgfplotstableread[skip first n=1]
{./numerical_tests_final/time_and_errors_projection_l2_2D_adm_2_degree_6_func_0_beta=0.005_threads_1_adaptive_refinement_wqMat_gaussRHS.dat}\wqdataDegSix

\begin{center}

\tikzset{every mark/.append style={scale=0.5}}

\begin{figure}[!h]
\begin{tikzpicture}[scale=0.5]
\begin{loglogaxis}[
    width=\textwidth,
    height=0.7\textwidth,
    title={2D case: DOFs vs. WQ-HR time details, degree $p=2$, admissibility $\admissibility=2$.},
    xlabel={$\ndofsHr$},
    ylabel={Time (sec.)},
    xmin=2e3,xmax=2.5e6,
    ymin=1e-1,ymax=3e2,
    legend pos=south east,
    grid=major
  ]
    
  \addplot[blue,mark=*,style=densely dotted]
    table [
      x index=2,
      y expr=\thisrowno{11}-\thisrowno{35}
      ] \wqdataDegTwo;
  \addlegendentry{Preprocessing}

  \addplot[red,mark=+]
    table [
      x index=2,
      y expr=\thisrowno{15}
      ] \wqdataDegTwo;
  \addlegendentry{Matrix formation}
  
  \addplot[black,mark=x,densely dashed]
    table [
      x index=2,
      y expr=\thisrowno{35}
      ] \wqdataDegTwo;
  \addlegendentry{Coeffs. evaluation}
\end{loglogaxis}
\end{tikzpicture}
\begin{tikzpicture}[scale=0.5]
\begin{loglogaxis}[
    width=\textwidth,
    height=0.7\textwidth,
    title={2D case: DOFs vs. WQ-HR time details, degree $p=3$, admissibility $\admissibility=2$.},
    xlabel={$\ndofsHr$},
    ylabel={Time (sec.)},
    xmin=2e3,xmax=2.5e6,
    ymin=1e-1,ymax=3e2,
    legend pos=south east,
    grid=major
  ]
    
  \addplot[blue,mark=*,style=densely dotted]
    table [
      x index=2,
      y expr=\thisrowno{11}-\thisrowno{35}
      ] \wqdataDegThree;
  \addlegendentry{Preprocessing}

  \addplot[red,mark=+]
    table [
      x index=2,
      y expr=\thisrowno{15}
      ] \wqdataDegThree;
  \addlegendentry{Matrix formation}
  
  \addplot[black,mark=x,densely dashed]
    table [
      x index=2,
      y expr=\thisrowno{35}
      ] \wqdataDegThree;
  \addlegendentry{Coeffs. evaluation}
\end{loglogaxis}
\end{tikzpicture}

\begin{tikzpicture}[scale=0.5]
\begin{loglogaxis}[
    width=\textwidth,
    height=0.7\textwidth,
    title={2D case: DOFs vs. WQ-HR time details, degree $p=4$, admissibility $\admissibility=2$.},
    xlabel={$\ndofsHr$},
    ylabel={Time (sec.)},
    xmin=2e3,xmax=2.5e6,
    ymin=1e-1,ymax=3e2,
    legend pos=south east,
    grid=major
  ]
    
  \addplot[blue,mark=*,style=densely dotted]
    table [
      x index=2,
      y expr=\thisrowno{11}-\thisrowno{35}
      ] \wqdataDegFour;
  \addlegendentry{Preprocessing}

  \addplot[red,mark=+]
    table [
      x index=2,
      y expr=\thisrowno{15}
      ] \wqdataDegFour;
  \addlegendentry{Matrix formation}
  
  \addplot[black,mark=x,densely dashed]
    table [
      x index=2,
      y expr=\thisrowno{35}
      ] \wqdataDegFour;
  \addlegendentry{Coeffs. evaluation}
\end{loglogaxis}
\end{tikzpicture}
\begin{tikzpicture}[scale=0.5]
\begin{loglogaxis}[
    width=\textwidth,
    height=0.7\textwidth,
    title={2D case: DOFs vs. WQ-HR time details, degree $p=5$, admissibility $\admissibility=2$.},
    xlabel={$\ndofsHr$},
    ylabel={Time (sec.)},
    xmin=2e3,xmax=2.5e6,
    ymin=1e-1,ymax=3e2,
    legend pos=south east,
    grid=major
  ]
    
  \addplot[blue,mark=*,style=densely dotted]
    table [
      x index=2,
      y expr=\thisrowno{11}-\thisrowno{35}
      ] \wqdataDegFive;
  \addlegendentry{Preprocessing}

  \addplot[red,mark=+]
    table [
      x index=2,
      y expr=\thisrowno{15}
      ] \wqdataDegFive;
  \addlegendentry{Matrix formation}
  
  \addplot[black,mark=x,densely dashed]
    table [
      x index=2,
      y expr=\thisrowno{35}
      ] \wqdataDegFive;
  \addlegendentry{Coeffs. evaluation}
\end{loglogaxis}
\end{tikzpicture}

\caption{$\ndofsHr$ vs.\ the CPU time needed to run the WQ-HR algorithm: preprocessing (Algorithm~\ref{algo:preprocessing}) and matrix computation (Algorithm~\ref{algo:assemble}) for the 2D case and admissibility parameter $\admissibility=2$. The cost for the matrix computation is split in the CPU cost for coefficient evaluations (\eqref{eq:coeffC} for $n=0,\dots,L$) and the CPU cost for executing the rest of the Algorithm~\ref{algo:assemble}.}
\label{fig:err_vs_time_WQ_details_2D_adm=2}
\end{figure}

\end{center}

%% file: numerical_tests_final/HB_WQ_tests_2D_adm3_dofs_vs_total_assembly_time.tex

\pgfplotstableread[skip first n=1]
{./numerical_tests_final/time_and_errors_projection_l2_2D_adm_3_degree_2_func_0_beta=0.005_threads_1_adaptive_refinement_wqMat_gaussRHS.dat}\wqdataDegTwo

\pgfplotstableread[skip first n=1]{./numerical_tests_final/time_and_errors_projection_l2_2D_adm_3_degree_2_func_0_beta=0.005_threads_1_adaptive_refinement_gaussMat_gaussRHS.dat}\gaussdataDegTwo

\pgfplotstableread[skip first n=1]
{./numerical_tests_final/time_and_errors_projection_l2_2D_adm_3_degree_3_func_0_beta=0.005_threads_1_adaptive_refinement_wqMat_gaussRHS.dat}\wqdataDegThree

\pgfplotstableread[skip first n=1]{./numerical_tests_final/time_and_errors_projection_l2_2D_adm_3_degree_3_func_0_beta=0.005_threads_1_adaptive_refinement_gaussMat_gaussRHS.dat}\gaussdataDegThree

\pgfplotstableread[skip first n=1]
{./numerical_tests_final/time_and_errors_projection_l2_2D_adm_3_degree_4_func_0_beta=0.005_threads_1_adaptive_refinement_wqMat_gaussRHS.dat}\wqdataDegFour

\pgfplotstableread[skip first n=1]{./numerical_tests_final/time_and_errors_projection_l2_2D_adm_3_degree_4_func_0_beta=0.005_threads_1_adaptive_refinement_gaussMat_gaussRHS.dat}\gaussdataDegFour

\pgfplotstableread[skip first n=1]
{./numerical_tests_final/time_and_errors_projection_l2_2D_adm_3_degree_5_func_0_beta=0.005_threads_1_adaptive_refinement_wqMat_gaussRHS.dat}\wqdataDegFive

\pgfplotstableread[skip first n=1]{./numerical_tests_final/time_and_errors_projection_l2_2D_adm_3_degree_5_func_0_beta=0.005_threads_1_adaptive_refinement_gaussMat_gaussRHS.dat}\gaussdataDegFive

\pgfplotstableread[skip first n=1]
{./numerical_tests_final/time_and_errors_projection_l2_2D_adm_3_degree_6_func_0_beta=0.005_threads_1_adaptive_refinement_wqMat_gaussRHS.dat}\wqdataDegSix

\pgfplotstableread[skip first n=1]{./numerical_tests_final/time_and_errors_projection_l2_2D_adm_3_degree_6_func_0_beta=0.005_threads_1_adaptive_refinement_gaussMat_gaussRHS.dat}\gaussdataDegSix

\begin{center}

\tikzset{every mark/.append style={scale=0.5}}

\begin{figure}[!h]
\begin{tikzpicture}
\begin{loglogaxis}[
    width=\textwidth,
    height=0.7\textwidth,
    title={2D case: DOFs vs. Total time for the matrix computation, admissibility $\admissibility=3$.},
    xlabel={$\ndofsHr$},
    ylabel={Time (sec.)},
    xmin=400,xmax=2.5e6,
    ymin=1e-1,ymax=1e3,
    legend pos=south east,
    legend columns=2, 
    legend style={
      /tikz/column 2/.style={column sep=15pt,},
    },
    grid=major
  ]
    


  \addplot[red,mark=*,mark size=1pt]
    table [
      x index=2,
      y expr=\thisrowno{11}+\thisrowno{15}
      ] \wqdataDegTwo;
   \addlegendentry{WQ-HR $p=2$}

  \addplot[red,mark=*,style=densely dotted]
    table [
      x index=2,
      y index=15
      ] \gaussdataDegTwo;
   \addlegendentry{Gauss $p=2$}



  \addplot[blue,mark=*,mark size=1pt]
    table [
      x index=2,
      y expr=\thisrowno{11}+\thisrowno{15}
      ] \wqdataDegThree;
   \addlegendentry{WQ-HR $p=3$}

  \addplot[blue,mark=*,style=densely dotted]
    table [
      x index=2,
      y index=15
      ] \gaussdataDegThree;
   \addlegendentry{Gauss $p=3$}



  \addplot[magenta,mark=*,mark size=1pt]
    table [
      x index=2,
      y expr=\thisrowno{11}+\thisrowno{15}
      ] \wqdataDegFour;
   \addlegendentry{WQ-HR $p=4$}

  \addplot[magenta,mark=*,style=densely dotted]
    table [
      x index=2,
      y index=15
      ] \gaussdataDegFour;
   \addlegendentry{Gauss $p=4$}




  \addplot[black,mark=*,mark size=1pt]
    table [
      x index=2,
      y expr=\thisrowno{11}+\thisrowno{15}
      ] \wqdataDegFive;
   \addlegendentry{WQ-HR $p=5$}

  \addplot[black,mark=*,style=densely dotted]
    table [
      x index=2,
      y index=15
      ] \gaussdataDegFive;
   \addlegendentry{Gauss $p=5$}



  \addplot[cyan,mark=*,mark size=1pt]
    table [
      x index=2,
      y expr=\thisrowno{11}+\thisrowno{15}
      ] \wqdataDegSix;
   \addlegendentry{WQ-HR $p=6$}

  \addplot[cyan,mark=*,style=densely dotted]
    table [
      x index=2,
      y index=15
      ] \gaussdataDegSix;
   \addlegendentry{Gauss $p=6$}


\end{loglogaxis}
\end{tikzpicture}
\caption{DOFs vs. total time for the matrix computation for the 2D case and admissibility parameter $\admissibility=3$. For WQ-HR the total time for the matrix computation is the sum of the time needed for the preprocessing (Section~\ref{sec:preprocessing}) and for the matrix formation (Section~\ref{sec:formation_algo}).}
\label{fig:dofs_vs_total_assembly_time_2D_adm=3}
\end{figure}
%
%

\end{center}

%% file: numerical_tests_final/HB_WQ_tests_3D_adm2_dofs_vs_total_assembly_time.tex

\pgfplotstableread[skip first n=1]{./numerical_tests_final/time_and_errors_projection_l2_3D_adm_2_degree_2_func_0_beta=0.1_threads_1_adaptive_refinement_wqMat_gaussRHS.dat}\wqdataDegTwo

\pgfplotstableread[skip first n=1]{./numerical_tests_final/time_and_errors_projection_l2_3D_adm_2_degree_2_func_0_beta=0.1_threads_1_adaptive_refinement_gaussMat_gaussRHS.dat}\gaussdataDegTwo

\pgfplotstableread[skip first n=1]{./numerical_tests_final/time_and_errors_projection_l2_3D_adm_2_degree_3_func_0_beta=0.1_threads_1_adaptive_refinement_wqMat_gaussRHS.dat}\wqdataDegThree

\pgfplotstableread[skip first n=1]{./numerical_tests_final/time_and_errors_projection_l2_3D_adm_2_degree_3_func_0_beta=0.1_threads_1_adaptive_refinement_gaussMat_gaussRHS.dat}\gaussdataDegThree

\pgfplotstableread[skip first n=1]{./numerical_tests_final/time_and_errors_projection_l2_3D_adm_2_degree_4_func_0_beta=0.1_threads_1_adaptive_refinement_wqMat_gaussRHS.dat}\wqdataDegFour

\pgfplotstableread[skip first n=1]{./numerical_tests_final/time_and_errors_projection_l2_3D_adm_2_degree_4_func_0_beta=0.1_threads_1_adaptive_refinement_gaussMat_gaussRHS.dat}\gaussdataDegFour

\pgfplotstableread[skip first n=1]{./numerical_tests_final/time_and_errors_projection_l2_3D_adm_2_degree_5_func_0_beta=0.1_threads_1_adaptive_refinement_wqMat_gaussRHS.dat}\wqdataDegFive

\pgfplotstableread[skip first n=1]{./numerical_tests_final/time_and_errors_projection_l2_3D_adm_2_degree_5_func_0_beta=0.1_threads_1_adaptive_refinement_gaussMat_gaussRHS.dat}\gaussdataDegFive



\begin{center}

\tikzset{every mark/.append style={scale=0.5}}

\begin{figure}[!h]
\begin{tikzpicture}
\begin{loglogaxis}[
    width=\textwidth,
    height=0.7\textwidth,
    title={3D case: DOFs vs. Total time for the matrix computation, admissibility $\admissibility=2$.},
    xlabel={$\ndofsHr$},
    ylabel={Time (sec.)},
    xmin=400,xmax=1e6,
    ymin=1e-1,ymax=2e3,
    legend pos=south east,
    legend columns=2, 
    legend style={
      /tikz/column 2/.style={column sep=15pt,},
    },
    grid=major
  ]
    
  \addplot[red,mark=*,mark size=1pt]
    table [
      x index=2,
      y expr=\thisrowno{11}+\thisrowno{15}
      ] \wqdataDegTwo;
  \addlegendentry{WQ-HR $p=2$}

  \addplot[red,mark=*,style=densely dotted]
    table [
      x index=2,
      y index=15
      ] \gaussdataDegTwo;
  \addlegendentry{Gauss $p=2$}

  \addplot[blue,mark=*,mark size=1pt]
    table [
      x index=2,
      y expr=\thisrowno{11}+\thisrowno{15}
      ] \wqdataDegThree;
  \addlegendentry{WQ-HR $p=3$}

  \addplot[blue,mark=*,style=densely dotted]
    table [
      x index=2,
      y index=15
      ] \gaussdataDegThree;
  \addlegendentry{Gauss $p=3$}

  \addplot[magenta,mark=*,mark size=1pt]
    table [
      x index=2,
      y expr=\thisrowno{11}+\thisrowno{15}
      ] \wqdataDegFour;
  \addlegendentry{WQ-HR $p=4$}

  \addplot[magenta,mark=*,style=densely dotted]
    table [
      x index=2,
      y index=15
      ] \gaussdataDegFour;
  \addlegendentry{Gauss $p=4$}

  \addplot[black,mark=*,mark size=1pt]
    table [
      x index=2,
      y expr=\thisrowno{11}+\thisrowno{15}
      ] \wqdataDegFive;
  \addlegendentry{WQ-HR $p=5$}

  \addplot[black,mark=*,style=densely dotted]
    table [
      x index=2,
      y index=15
      ] \gaussdataDegFive;
  \addlegendentry{Gauss $p=5$}



\end{loglogaxis}
\end{tikzpicture}
\caption{DOFs vs. total time for the matrix computation for the 3D case and admissibility parameter $\admissibility=2$. For WQ-HR the total time for the matrix computation is the sum of the time needed for the preprocessing (Section~\ref{sec:preprocessing}) and for the matrix formation (Section~\ref{sec:formation_algo}).}
\label{fig:dofs_vs_total_assembly_time_3D_adm=2}
\end{figure}
%
%

\end{center}

%% file: numerical_tests_final/HB_WQ_tests_3D_adm2_L2error_vs_total_assembly_time.tex

\pgfplotstableread[skip first n=1]{./numerical_tests_final/time_and_errors_projection_l2_3D_adm_2_degree_2_func_0_beta=0.1_threads_1_adaptive_refinement_wqMat_gaussRHS.dat}\wqdataDegTwo

\pgfplotstableread[skip first n=1]{./numerical_tests_final/time_and_errors_projection_l2_3D_adm_2_degree_2_func_0_beta=0.1_threads_1_adaptive_refinement_gaussMat_gaussRHS.dat}\gaussdataDegTwo

\pgfplotstableread[skip first n=1]{./numerical_tests_final/time_and_errors_projection_l2_3D_adm_2_degree_3_func_0_beta=0.1_threads_1_adaptive_refinement_wqMat_gaussRHS.dat}\wqdataDegThree

\pgfplotstableread[skip first n=1]{./numerical_tests_final/time_and_errors_projection_l2_3D_adm_2_degree_3_func_0_beta=0.1_threads_1_adaptive_refinement_gaussMat_gaussRHS.dat}\gaussdataDegThree

\pgfplotstableread[skip first n=1]{./numerical_tests_final/time_and_errors_projection_l2_3D_adm_2_degree_4_func_0_beta=0.1_threads_1_adaptive_refinement_wqMat_gaussRHS.dat}\wqdataDegFour

\pgfplotstableread[skip first n=1]{./numerical_tests_final/time_and_errors_projection_l2_3D_adm_2_degree_4_func_0_beta=0.1_threads_1_adaptive_refinement_gaussMat_gaussRHS.dat}\gaussdataDegFour

\pgfplotstableread[skip first n=1]{./numerical_tests_final/time_and_errors_projection_l2_3D_adm_2_degree_5_func_0_beta=0.1_threads_1_adaptive_refinement_wqMat_gaussRHS.dat}\wqdataDegFive

\pgfplotstableread[skip first n=1]{./numerical_tests_final/time_and_errors_projection_l2_3D_adm_2_degree_5_func_0_beta=0.1_threads_1_adaptive_refinement_gaussMat_gaussRHS.dat}\gaussdataDegFive



\begin{center}

\tikzset{every mark/.append style={scale=0.5}}

\begin{figure}[!h]
\begin{tikzpicture}
\begin{loglogaxis}[
    width=\textwidth,
    height=0.7\textwidth,
    title={3D case: $L^{2}$-error vs. Total time for the matrix computation, admissibility $\admissibility=2$.},
    xlabel={$L^{2}$-error},
    ylabel={Time (sec.)},
    xmin=9e-7,xmax=1.5e-1,
    ymin=1e-1,ymax=2e3,
    legend pos=south west,
    legend columns=2, 
    legend style={
      /tikz/column 2/.style={column sep=15pt,},
    },
    grid=major
  ]
    
  \addplot[red,mark=*,mark size=1pt]
    table [
      x index=4,
      y expr=\thisrowno{11}+\thisrowno{15}
      ] \wqdataDegTwo;
  \addlegendentry{WQ-HR $p=2$}

  \addplot[red,mark=*,style=densely dotted]
    table [
      x index=4,
      y index=15
      ] \gaussdataDegTwo;
  \addlegendentry{Gauss $p=2$}   

  \addplot[blue,mark=*,mark size=1pt]
    table [
      x index=4,
      y expr=\thisrowno{11}+\thisrowno{15}
      ] \wqdataDegThree;
  \addlegendentry{WQ-HR $p=3$}

  \addplot[blue,mark=*,style=densely dotted]
    table [
      x index=4,
      y index=15
      ] \gaussdataDegThree;
  \addlegendentry{Gauss $p=3$}

  \addplot[magenta,mark=*,mark size=1pt]
    table [
      x index=4,
      y expr=\thisrowno{11}+\thisrowno{15}
      ] \wqdataDegFour;
  \addlegendentry{WQ-HR $p=4$}

  \addplot[magenta,mark=*,style=densely dotted]
    table [
      x index=4,
      y index=15
      ] \gaussdataDegFour;
  \addlegendentry{Gauss $p=4$}

  \addplot[black,mark=*,mark size=1pt]
    table [
      x index=4,
      y expr=\thisrowno{11}+\thisrowno{15}
      ] \wqdataDegFive;
  \addlegendentry{WQ-HR $p=5$}

  \addplot[black,mark=*,style=densely dotted]
    table [
      x index=4,
      y index=15
      ] \gaussdataDegFive;
  \addlegendentry{Gauss $p=5$}



\end{loglogaxis}
\end{tikzpicture}
\caption{$L^{2}$-error vs. total time for the matrix computation for the 3D case and admissibility parameter $\admissibility=2$. For WQ-HR the total time for the matrix computation is the sum of the time needed for the preprocessing (Section~\ref{sec:preprocessing}) and for the matrix formation (Section~\ref{sec:formation_algo}).}
\label{fig:err_vs_total_assembly_time_3D_adm=2}
\end{figure}
%
%

\end{center}

%% file: numerical_tests_final/HB_WQ_tests_3D_adm2_dofs_vs_WQ_details.tex

\pgfplotstableread[skip first n=1]{./numerical_tests_final/time_and_errors_projection_l2_3D_adm_2_degree_2_func_0_beta=0.1_threads_1_adaptive_refinement_wqMat_gaussRHS.dat}\wqdataDegTwo

\pgfplotstableread[skip first n=1]{./numerical_tests_final/time_and_errors_projection_l2_3D_adm_2_degree_3_func_0_beta=0.1_threads_1_adaptive_refinement_wqMat_gaussRHS.dat}\wqdataDegThree

\pgfplotstableread[skip first n=1]{./numerical_tests_final/time_and_errors_projection_l2_3D_adm_2_degree_4_func_0_beta=0.1_threads_1_adaptive_refinement_wqMat_gaussRHS.dat}\wqdataDegFour

\pgfplotstableread[skip first n=1]{./numerical_tests_final/time_and_errors_projection_l2_3D_adm_2_degree_5_func_0_beta=0.1_threads_1_adaptive_refinement_wqMat_gaussRHS.dat}\wqdataDegFive


\begin{center}

\tikzset{every mark/.append style={scale=0.5}}

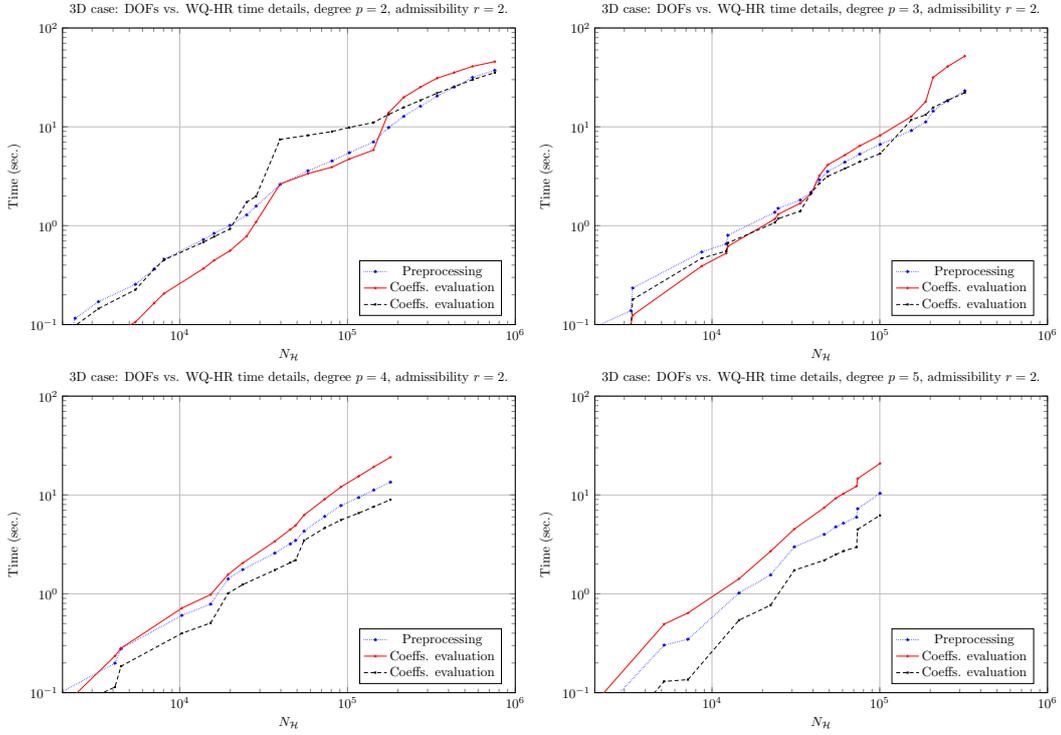
\begin{figure}[!h]
\begin{tikzpicture}[scale=0.5]
\begin{loglogaxis}[
    width=\textwidth,
    height=0.7\textwidth,
    title={3D case: DOFs vs. WQ-HR time details, degree $p=2$, admissibility $\admissibility=2$.},
    xlabel={$\ndofsHr$},
    ylabel={Time (sec.)},
    xmin=2e3,xmax=1e6,
    ymin=1e-1,ymax=1e2,
    legend pos=south east,
    grid=major
  ]
    
  \addplot[blue,mark=*,style=densely dotted]
    table [
      x index=2,
      y expr=\thisrowno{11}-\thisrowno{35}
      ] \wqdataDegTwo;
  \addlegendentry{Preprocessing}

  \addplot[red,mark=+]
    table [
      x index=2,
      y expr=\thisrowno{15}
      ] \wqdataDegTwo;
  \addlegendentry{Coeffs. evaluation}
  
  \addplot[black,mark=x,densely dashed]
    table [
      x index=2,
      y expr=\thisrowno{35}
      ] \wqdataDegTwo;
  \addlegendentry{Coeffs. evaluation}
\end{loglogaxis}
\end{tikzpicture}
\begin{tikzpicture}[scale=0.5]
\begin{loglogaxis}[
    width=\textwidth,
    height=0.7\textwidth,
    title={3D case: DOFs vs. WQ-HR time details, degree $p=3$, admissibility $\admissibility=2$.},
    xlabel={$\ndofsHr$},
    ylabel={Time (sec.)},
    xmin=2e3,xmax=1e6,
    ymin=1e-1,ymax=1e2,
    legend pos=south east,
    grid=major
  ]
    
  \addplot[blue,mark=*,style=densely dotted]
    table [
      x index=2,
      y expr=\thisrowno{11}-\thisrowno{35}
      ] \wqdataDegThree;
  \addlegendentry{Preprocessing}

  \addplot[red,mark=+]
    table [
      x index=2,
      y expr=\thisrowno{15}
      ] \wqdataDegThree;
  \addlegendentry{Coeffs. evaluation}
  
  \addplot[black,mark=x,densely dashed]
    table [
      x index=2,
      y expr=\thisrowno{35}
      ] \wqdataDegThree;
  \addlegendentry{Coeffs. evaluation}
\end{loglogaxis}
\end{tikzpicture}

\begin{tikzpicture}[scale=0.5]
\begin{loglogaxis}[
    width=\textwidth,
    height=0.7\textwidth,
    title={3D case: DOFs vs. WQ-HR time details, degree $p=4$, admissibility $\admissibility=2$.},
    xlabel={$\ndofsHr$},
    ylabel={Time (sec.)},
    xmin=2e3,xmax=1e6,
    ymin=1e-1,ymax=1e2,
    legend pos=south east,
    grid=major
  ]
    
  \addplot[blue,mark=*,style=densely dotted]
    table [
      x index=2,
      y expr=\thisrowno{11}-\thisrowno{35}
      ] \wqdataDegFour;
  \addlegendentry{Preprocessing}

  \addplot[red,mark=+]
    table [
      x index=2,
      y expr=\thisrowno{15}
      ] \wqdataDegFour;
  \addlegendentry{Coeffs. evaluation}
  
  \addplot[black,mark=x,densely dashed]
    table [
      x index=2,
      y expr=\thisrowno{35}
      ] \wqdataDegFour;
  \addlegendentry{Coeffs. evaluation}
\end{loglogaxis}
\end{tikzpicture}
\begin{tikzpicture}[scale=0.5]
\begin{loglogaxis}[
    width=\textwidth,
    height=0.7\textwidth,
    title={3D case: DOFs vs. WQ-HR time details, degree $p=5$, admissibility $\admissibility=2$.},
    xlabel={$\ndofsHr$},
    ylabel={Time (sec.)},
    xmin=2e3,xmax=1e6,
    ymin=1e-1,ymax=1e2,
    legend pos=south east,
    grid=major
  ]
    
  \addplot[blue,mark=*,style=densely dotted]
    table [
      x index=2,
      y expr=\thisrowno{11}-\thisrowno{35}
      ] \wqdataDegFive;
  \addlegendentry{Preprocessing}

  \addplot[red,mark=+]
    table [
      x index=2,
      y expr=\thisrowno{15}
      ] \wqdataDegFive;
  \addlegendentry{Coeffs. evaluation}
  
  \addplot[black,mark=x,densely dashed]
    table [
      x index=2,
      y expr=\thisrowno{35}
      ] \wqdataDegFive;
  \addlegendentry{Coeffs. evaluation}
\end{loglogaxis}
\end{tikzpicture}

\caption{DOFs vs. the CPU time needed to run the WQ-HR algorithm: preprocessing (Algorithm~\ref{algo:preprocessing}) and matrix computation (Algorithm~\ref{algo:assemble}) for the 3D case and admissibility parameter $\admissibility=2$. The cost for the matrix computation is split in the CPU cost for coefficient evaluations (Eq.~\eqref{eq:coeffC}, for $n=0,\dots,L$) and the CPU cost for executing the rest of the Algorithm~\ref{algo:assemble}.}
\label{fig:err_vs_time_WQ_details_3D_adm=2}
\end{figure}

\end{center}